\documentclass[a4paper,12pt]{amsart}
\usepackage{amssymb}
\usepackage{ifthen}
\usepackage{graphicx}
\usepackage{color,xcolor}
\usepackage{mathrsfs}
\numberwithin{equation}{section}
\newtheorem{thm}{Theorem}[section]
\newtheorem{lem}{Lemma}[section]
\newtheorem{cor}{Corollary}[section]
\newtheorem{prop}{Proposition}[section]
\newtheorem{defn}{Definition}[section]
\newtheorem{example}{Example}%[section]

\newtheorem{claim}{Claim}

\newtheorem{conj}{Conjecture}
\newtheorem{prob}{Problem}

\newenvironment{rem}{%
\bigskip
\noindent \textsl{{\sl Remark. }}}{\bigskip}
\newenvironment{rems}{%
\bigskip
\noindent \textsl{{\sl Remarks. }}}{\bigskip}

\newenvironment{pf}[1][]{%
 \vskip 1mm
 \noindent
 \ifthenelse{\equal{#1}{}}%
  {{\slshape Proof. }}%
  {{\slshape #1.} }%
 }%
{\qed\medskip}

\newcounter{alphabet}
\newcounter{tmp}

\makeatletter
\newcommand{\Ref}[1]{\@ifundefined{r@#1}{}{\setcounter{tmp}{\ref{#1}}\Alph{tmp}}}
\makeatother
\DeclareMathOperator{\RE}{Re}

\newcommand{\ID}{{\mathbb D}}

\def\be{\begin{equation}}
\def\ee{\end{equation}}
\def\bes{\begin{equation*}}
\def\ees{\end{equation*}}

\newcommand{\bee}{\begin{enumerate}}
\newcommand{\eee}{\end{enumerate}}

\newcommand{\blem}{\begin{lem}}
\newcommand{\elem}{\end{lem}}
\newcommand{\bthm}{\begin{thm}}
\newcommand{\ethm}{\end{thm}}
\newcommand{\bcor}{\begin{cor}}
\newcommand{\ecor}{\end{cor}}
\newcommand{\beg}{\begin{example}}
\newcommand{\eeg}{\end{example}}
\newcommand{\begs}{\begin{examples}}
\newcommand{\eegs}{\end{examples}}
\newcommand{\bdefe}{\begin{defn}}
\newcommand{\edefe}{\end{defn}}
\newcommand{\bprob}{\begin{prob}}
\newcommand{\eprob}{\end{prob}}
\newcommand{\bques}{\begin{ques}}
\newcommand{\eques}{\end{ques}}
\newcommand{\bei}{\begin{itemize}}
\newcommand{\eei}{\end{itemize}}

\newcommand{\bde}{\begin{deter}}
\newcommand{\ede}{\end{deter}}
\newcommand{\bca}{\begin{case}}
\newcommand{\eca}{\end{case}}
\newcommand{\bcl}{\begin{claim}}
\newcommand{\ecl}{\end{claim}}
\newcommand{\bcon}{\begin{conj}}
\newcommand{\econ}{\end{conj}}
\newcommand{\bcons}{\begin{conjs}}
\newcommand{\econs}{\end{conjs}}
\newcommand{\bprop}{\begin{propo}}
\newcommand{\eprop}{\end{propo}}
\newcommand{\br}{\begin{rem}}
\newcommand{\er}{\end{rem}}
\newcommand{\brs}{\begin{rems}}
\newcommand{\ers}{\end{rems}}
\newcommand{\bo}{\begin{obser}}
\newcommand{\eo}{\end{obser}}
\newcommand{\bos}{\begin{obsers}}
\newcommand{\eos}{\end{obsers}}
\newcommand{\bpf}{\begin{pf}}
\newcommand{\epf}{\end{pf}}
\newcommand{\ba}{\begin{array}}
\newcommand{\ea}{\end{array}}
\newcommand{\beq}{\begin{eqnarray}}
\newcommand{\beqq}{\begin{eqnarray*}}
\newcommand{\eeq}{\end{eqnarray}}
\newcommand{\eeqq}{\end{eqnarray*}}

\newcommand{\ds}{\displaystyle}

\newcounter{minutes}
\divide\time by 60
\newcounter{hours}
\multiply\time by 60 \addtocounter{minutes}{-\time}
\begin{document}
\title[Uniformly locally univalent  harmonic mappings associated with the pre-Schwarzian norm]
{Uniformly locally univalent harmonic mappings associated with the pre-Schwarzian norm}

%=========================================================================
\thanks{%$^\dagger$
File:~\jobname .tex,
          printed: \number\day-\number\month-\number\year,
          \thehours.\ifnum\theminutes<10{0}\fi\theminutes}
%=========================================================================

\author[G. Liu]{Gang Liu}
\address{G. Liu, College of Mathematics and Statistics
 (Hunan Provincial Key Laboratory of Intelligent Information Processing and Application),
Hengyang Normal University, Hengyang,  Hunan 421002, People's Republic of China.}
\email{liugangmath@sina.cn}

\author[S. Ponnusamy]{Saminathan Ponnusamy
%$^\dagger $
%${}^{~\mathbf{*}}$
}
%\address{Saminathan  Ponnusamy,
%Indian Statistical Institute (ISI), Chennai Centre, SETS (Society
%for Electronic Transactions and Security), MGR Knowledge City, CIT
%Campus, Taramani, Chennai 600 113, India.}
%
\address{S. Ponnusamy, Department of Mathematics,
Indian Institute of Technology Madras, Chennai-600 036, India}
\email{samy@iitm.ac.in, samy@isichennai.res.in }

%\address{S. Ponnusamy, Stat-Math Unit,
%Indian Statistical Institute (ISI), Chennai Centre,
%110, Nelson Manickam Road,
%Aminjikarai, Chennai, 600 029, India.}
%\email{samy@isichennai.res.in, samy@iitm.ac.in}

\subjclass[2010]{Primary: 30H10, 30H30, 30H35, 31A05; Secondary 30C55}
%; Secondary: 30C20, 30C55}
%,  31B05, 31C05}

\keywords{Uniformly locally univalent, analytic functions, harmonic and stable harmonic mappings, pre-Schwarzian and Schwarzian,
uniformly hyperbolic radius, distortion, Hardy space, Bloch space, BMOA, BMOH, subordination.
%\\
%$%{}^{\mathbf{*}}
%^\dagger$ {\tt This author is on leave from the Department of Mathematics,
%Indian Institute of Technology Madras, Chennai-600 036, India}
}
%\thanks{ }
%\maketitle

%\newfont{\Bbb}{msbm10 scaled\magstephalf}
\begin{abstract}
In this paper, we consider the class of uniformly locally univalent harmonic mappings in the unit disk
and build a relationship between its pre-Schwarzian norm and uniformly hyperbolic radius. Also, we establish
eight ways of characterizing uniformly locally univalent sense-preserving harmonic mappings.
We also present some sharp distortions and growth estimates and investigate their connections with Hardy spaces.
Finally, we study subordination principles of norm estimates.
\end{abstract}

\maketitle \pagestyle{myheadings}
\markboth{G. Liu and S. Ponnusamy}{Harmonic mappings, the pre-Schwarzian and the Bloch space}

\section{Introduction}  \label{LS3-sec1}

Harmonic mappings play an important role in various branches of applied mathematics including the study of liquid crystals,
both in theory and in practice. There are many classical
approaches to deal with harmonic maps in various settings. For example, A. Aleman and A. Constantin \cite{AleConst2012}
developed tools using complex analytic theory and the univalence of the labelling map to solve fluid flow problems in
a surprisingly simple form.  More recently, O. Constantin and M. J. Mart\'{i}n
\cite{ConstMartin2017} proposed a new approach to obtain a
complete solution to the problem of classifying all two dimensional ideal fluid flows with harmonic labelling maps.
This approach is based on ideas from the theory of harmonic mappings by finding two harmonic maps with same Jacobians
and illustrates the deep links between the fields of complex analysis and fluid mechanics.
Investigations of this type have prompted renewed interest in the study of sense-preserving harmonic mappings.
The present article is concerned with Schwarzian and pre-Schwarzian norms  defined in the unit disk, and in particular, with
certain important function spaces. In addition, we introduce several new ideas and tools for a number of problems
in the case of harmonic mappings.

\subsection{Basic notations}
%All the results obtained in the paper regarding sense-preserving harmonic mappings can be adapted to sense-reversing ones.
 A complex-valued function $f$ in the unit disk $\ID=\{z:\,|z|<1\}$ is called a {\it harmonic mapping}
if it satisfies the Laplace equation %$\Delta f=0$,
%where $\Delta$ denotes  the complex Laplacian operator
$\Delta f=4f_{z\overline{z}}=f_{xx}+f_{yy}=0$.
It is known that $f$ has  a canonical representation $f=h+\overline{g}$ with $g(0)=0$,
where $h$ and $g$ are analytic functions in $\mathbb{D}$ and $J_f=|h'|^2-|g'|^2$ denotes the Jacobian of $f$.
As is usual, we call $h$ the {\it analytic part} of $f$
and $g$ the {\it co-analytic part} of $f$.
Lewy \cite{lew} proved that $f=h+\overline{g}$ is locally univalent in $\mathbb{D}$
if and only if $J_f(z)\neq 0$ in $\ID$.
Without loss of generality, we consider harmonic mappings $f$ that are sense-preserving, i.e. $J_f>0$ or
equivalently $|h'|>|g'|$ in $\mathbb{D}$.
%In view of this, a locally univalent harmonic mapping $f$ is called {\it sense-preserving}
%and {\it sense-reversing}  in $\mathbb{D}$ according to $J_f>0$ and $J_f<0$ in $\mathbb{D}$, respectively.
%Since $J_f=-J_{\overline{f}}$, $f$ is sense-preserving if and only if $\overline{f}$ is sense-reversing.
%Thus, $f=h+\overline{g}$ is sense-preserving in $\mathbb{D}$ if and only if $|h'|>|g'|$ in $\mathbb{D}$.
In this case,  its dilatation $\omega_f=g'/h'$ has the property that $|\omega_f|<1$ in $\mathbb{D}$.
Especially, if $|\omega_f|\leq k<1$ in $\mathbb{D}$, then $f$ is called a $K-$quasiconformal mapping, where $K=(1+k)/(1-k)$.
More details about planar harmonic mappings, may be found in the monograph of Duren \cite{dur2004} and in %the
the survey article of Ponnusamy and Rasila \cite{PR}.

%Throughout this paper we consider sense-preserving complex-valued harmonic mappings defined on the unit disk $\mathbb{D}=\{z:\,|z|<1\}$.
For the convenience of the reader, we first list down the following notations and terminologies whose precise definitions
will be presented at appropriate places.
\bei
\item ULU (ULC) - uniformly locally univalent (uniformly locally convex)
\item SAULU - stable analytic uniformly locally univalent
\item SAULC - stable analytic uniformly locally convex
\item SHULU  - stable harmonic uniformly locally univalent
\item SHULC - stable harmonic uniformly locally convex
\item SHU (SHC) - stable harmonic univalent (stable harmonic convex)
\item PSD (SD) - pre-Schwarzian derivative  (Schwarzian derivative)
\item PSN (SN) - pre-Schwarzian norm (Schwarzian norm)
\item SBAPSN  - stable bounded analytic pre-Schwarzian norm
\item SBASN - stable bounded analytic Schwarzian norm
\item SBHPSN  - stable bounded harmonic pre-Schwarzian norm
\item SBHSN - stable bounded harmonic Schwarzian norm
%\item $\mathcal{A}=\{f:\, \mbox{$f$ is analytic in $\mathbb{D}$ with the normalizations $f(0)=f'(0)-1=0$} $
\item $\mathcal{H}=\left \{f=h+\overline{g}:\, \mbox{$f$ is a sense-preserving harmonic mapping in $\mathbb{D}$ satisfying }\right .$\\
~~~~~ $\left . ~~~\mbox{  the normalizations $h(0)=h'(0)-1=g(0)=0$}\right \}$
\item $\mathcal{H}_0=\{f=h+\overline{g}\in\mathcal{H}:\, g'(0)=0\}$
\eei
Sometimes we write $f\in {\rm ULU}$ to convey that $f$ is a uniformly locally univalent function in $\ID$. Similar convention will be followed for
other cases.
%Throughout this paper we consider sense-preserving complex-valued harmonic mappings defined on the unit disk $\mathbb{D}=\{z:\,|z|<1\}$.

\subsection{ULU harmonic mappings}  \label{LS3-sec1.1}
%Next we introduce the notation of %{\it uniformly locally univalent}
%ULU harmonic mappings in $\mathbb{D}$.
%We denote $d_h(z,a)$ as the hyperbolic distance of $z,~a\in\mathbb{D}$,
%that is
Let $z,~a\in\mathbb{D}$. We denote the hyperbolic distance between $z$ and $a$ by
$$d_h(z,a)=2\tanh^{-1}\left|\frac{z-a}{1-\overline{a}z}\right|.
$$
The hyperbolic disk in $\mathbb{D}$ with center $a\in\mathbb{D}$
and  hyperbolic radius $\rho$, $0<\rho\leq\infty$, is defined by
$$D_h(a,\rho)=\{z\in\mathbb{D}:\,d_h(z,a)<\rho\}.
$$
We say that a sense-preserving harmonic mapping $f=h+\overline{g}$ in $\mathbb{D}$ is a ULU harmonic mapping in $\mathbb{D}$ if $\rho(f)>0$,
where
$$\rho(f)=\inf_{z\in\mathbb{D}}\left\{\sup_{\rho_z>0}\{\rho_z:\,f ~\mbox{is univalent in}~ D_h(z,\rho_z)\}\right\}.
$$
%, then we say that $f$ is  ULU in $\mathbb{D}$
The number $\rho(f)$ is called the  {\it uniformly hyperbolic radius} of $f$. Moreover,
$f$ is univalent in $\mathbb{D}$ if and only if $\rho(f)=\infty$.

\subsection{PSD and PSN of harmonic mappings}  \label{LS3-sec1.2}

Let $f=h+\overline{g}$ be a sense-preserving harmonic mapping in $\mathbb{D}$ with $\omega:= \omega _f=g'/h'$.
%The {\it pre-Schwarzian derivative and norm} of $f$ are defined by
Then the PSD and the PSN of $f$ are defined by
\begin{equation}\label{LS3-equ1.1}
P_f=(\log J_f)_z=\frac{h''\overline{h'}-g''\overline{g'}}{|h'|^2-|g'|^2}=\frac{h''}{h'}-\frac{\overline{\omega}\omega'}{1-|\omega|^2}
 %\quad \mbox{for} \quad z\in\mathbb{D},
\end{equation}
and
\begin{equation*}
||P_{f}||=\sup_{z\in\mathbb{D}}(1-|z|^2)|P_{f}(z)|,
\end{equation*}
respectively. Clearly, in the analytic case $\omega $ in \eqref{LS3-equ1.1}
has taken to be identically $0$ in $\ID$, and thus, throughout we
use the same notations for the PSD and the PSN in the case of analytic functions as well.

%Since $J_{\overline{f}}=-J_f$, we have $P_{\overline{f}}=P_f$ and $||P_{\overline{f}}||=||P_{f}||$,
%and thus it suffices to consider  sense-preserving harmonic mappings. Consequently,
The PSD has affine invariance property:
\begin{equation}\label{LS3-equ1.2}
\quad P_f=P_{A\circ f}, \quad
A(z)=az+\overline{bz}+c, \quad a,~b,~c\in\mathbb{C}\quad \mbox{and}\quad |a|>|b|.
\end{equation}
% in $\mathbb{D}$, where $A(z)=az+\overline{bz}+c$, $a,~b,~c\in\mathcal{C}$ and $|a|>|b|$.
Note that $A\circ f$ is still sense-preserving in $\mathbb{D}$.
%Furthermore, when $f$ is sense-preserving and $|a|>|b|$, $A\circ f$ is still sense-preserving in $\mathbb{D}$
%because the dilatation
%\begin{equation}\label{LS3-equ1.5}
%\omega_{A\circ f}=\frac{b+\overline{a}\omega_f}{a+\overline{b}\omega_f}
%\end{equation}
%has the property $|\omega_{A\circ f}|<1$ in $\mathbb{D}$.

The  above definitions of the PSD and the PSN for harmonic mappings were introduced by Hern\'{a}ndez and Mart\'{\i}n in \cite{HM2015} (see also \cite{CDO}),
which coincide with the corresponding analytic definitions (see \cite{dur1983,pom1992}). It is well known that the PSN of a locally univalent analytic function
is an important quantity  in the study of the global univalence. For example, if $f$ is a univalent analytic function in  $\mathbb{D}$,
then  $||P_{f}||\leq6$, which is sharp.
%The well-known inequality is consequence of the classical result
%by Bieberbach for univalent analytic function $f\in\mathcal{A}$.
Conversely, if $||P_{f}||\leq1$ holds for a locally univalent analytic function $f$ in $\mathbb{D}$,
then $f$ is necessarily univalent in $\ID$ and the constant $1$ is sharp (see \cite{bec,BP}).
Recently, new criteria for the univalence
of  harmonic mappings in terms of the PSD or the PSN have been established in \cite{ANS, gra, HM2015}.

%\subsection{Relationship between the ULU and the PSN}  \label{LS3-sec1.3}
\subsection{Relationship between ULU and PSN}  \label{LS3-sec1.3}
%If the ULU harmonic mapping $f$ in $\mathbb{D}$  is restricted to be an analytic function,
%then it reduces to the ULU analytic function in $\mathbb{D}$.
%Such function has a connection with its pre-Schwarzian  norm adhered to its pre-Schwarzian derivative.
%To avoid confusion of symbols, let $h$ be a locally univalent analytic function in $\mathbb{D}.$
%The {\it pre-Schwarzian  derivative and norm of $h$} are defined by
%\begin{equation}\label{LS3-equ1.1}
% P_h(z)=\frac{h''}{h'},\quad z\in\mathbb{D}, \quad \text{and}
%  \quad ||P_{h}||=\sup_{z\in\mathbb{D}}|P_{h}(z)|(1-|z|^2),
%   \end{equation}
%respectively.
Yamashita  \cite{yam1977} showed that a locally univalent analytic function $f$ in $\mathbb{D}$ is ULU in $\mathbb{D}$
if and only if $||P_{f}||$ is bounded.
Later, Kim and Sugawa \cite{KS2002} investigated the growth of various quantities
for a ULU analytic function in $\mathbb{D}$ by means of finite the PSN.
Since $P_{\phi\circ f}=P_f$ for any linear transformation $\phi(z)=az+b~(a\neq0)$,
they just considered the following normalized function space
$$\mathcal{B}_{A}=\{f\in\mathcal{A}:\,||P_f||<\infty\},
$$
where $\mathcal{A}$ is the set of analytic functions $f$ in $\mathbb{D}$ with the normalizations $f(0)=f'(0)-1=0$.
In fact, the space $\mathcal{B}_{A}$ has the structure of a nonseparable
complex Banach space under the Hornich operation (see \cite{yam1975}).
To obtain some precise results, it was necessary to study the subset of
$\mathcal{B}_{A}$:
$$\mathcal{B}_{A}(\lambda)=\{f\in\mathcal{A}:\,||P_f||\leq2\lambda\},
$$
where $\lambda\geq0$ and the factor 2 is due to only some technical reason.

Following the proof of \cite[Theorem~7]{HM2015}, we see that a sense-preserving harmonic mapping $f$ in $\mathbb{D}$
is ULU in $\mathbb{D}$ if and only if $||P_{f}||$ is bounded, which will be also proved in Section \ref{LS3-sec3} by other method.
Therefore, the primary aim of this paper is to extend some of the results from \cite{KS2002}
to sense-preserving and ULU harmonic mappings in $\mathbb{D}$ associated with finite the PSN.
Since the PSD preserves affine invariance, in what follows,
we only to consider the following set of normalized functions:
$$\mathcal{B}_H=\{f\in\mathcal{H}:\, ||P_f||<\infty\}.
$$
%where $\mathcal{H}$ is the set of normalized sense-preserving harmonic mappings $f=h+\overline{g}$
%in $\mathbb{D}$ satisfying $h(0)=h'(0)-1=g(0)=0$.
%Let $\mathcal{H}_0$ be the subset of  $\mathcal{H}$ with $g'(0)=0$.
If we concern only on the PSN, then
$\mathcal{B}_H$ can be further restricted to be $\mathcal{B}_H^0:=\mathcal{B}_H\cap\mathcal{H}_0$.
In fact, if $f=h+\overline{g}\in\mathcal{B}_H$ and $A(z)=\frac{z-\overline{b_1z}}{1-|b_1|^2}$ $(b_1=g'(0))$,
then it follows from \eqref{LS3-equ1.2} that $||P_{A\circ f}||=||P_f||$
and it is also easy to see that $A\circ f\in\mathcal{B}_H^0$.

Let $f=h+\overline{g}$ be a sense-preserving harmonic mapping in $\mathbb{D}$.
Then, motivated by the works of \cite{HM2013},
in Section \ref{LS3-sec2}, we will build some sharp inequalities between  $||P_{h+\varepsilon_1\overline{g}}||$
and $||P_{h+\varepsilon_2 g}||$, where $\varepsilon_1,\varepsilon_2\in\overline{\mathbb{D}}$.
In particular, we obtain the following important implication:
\begin{equation}\label{LS3-equ1.3}
f\in\mathcal{B}_H(\lambda):=\{f\in\mathcal{H}:\, ||P_f||\leq\lambda\}\Rightarrow \frac{ h+\varepsilon g}{1+\varepsilon g'(0)}
\in\mathcal{B}_{A}\left(\frac{\lambda+1}{2}\right)
\quad \forall~\varepsilon\in\overline{\mathbb{D}},
\end{equation}
%Here,  and
%$$\mathcal{B}_H(\lambda)=\{f\in\mathcal{H}:\, ||P_f||\leq\lambda\}.
%$$
where $\lambda\geq0$.
In Section \ref{LS3-sec3}, for any given sense-preserving and ULU harmonic mapping in  the unit disk,
we give a relationship between its PSN and uniformly hyperbolic radius.
Combining the above results with some works about ULU harmonic mappings,
plenty of equivalent conditions for a sense-preserving and ULU harmonic mapping
in  the unit disk are obtained in Section \ref{LS3-sec4}.
To present some  sharp examples in Sections \ref{LS3-sec6} and \ref{LS3-sec7}, we introduce a class of
sense-preserving harmonic mappings with prescribed PSN in Section \ref{LS3-sec5}.
These results help us to obtain sharp distortion,
growth and covering theorems for $\mathcal{B}_H(\lambda)$
or $\mathcal{B}_H^0(\lambda):=\mathcal{B}_H(\lambda)\cap\mathcal{H}_0$  in Section \ref{LS3-sec6}.
Applying \eqref{LS3-equ1.3} and the corresponding results in \cite{KS2002} and \cite{PQW},
the growth of coefficients and the relationship with Hardy space
for the class $\mathcal{B}_H(\lambda)$ are considered in Sections \ref{LS3-sec7} and \ref{LS3-sec8}, respectively.
Finally, some subordination principles of the PSN estimates are also obtained in Section \ref{LS3-sec9}.

\section{Some inequalities Concerning Pre-Schwarzian norm} \label{LS3-sec2}

We now state our key inequalities which will provide important connections
between ULU analytic functions and ULU harmonic mappings in the unit disk.

\bthm\label{LS3-thm2.1}
Let $f=h+\overline{g}$ be a sense-preserving harmonic mapping in $\mathbb{D}$.
Then either $||P_{h+\varepsilon g}||=||P_f||=\infty$ or both $||P_{h+\varepsilon g}||$ and $||P_f||$
are finite for each $\varepsilon\in\overline{\mathbb{D}}$.
If $||P_f||<\infty$, then the inequality
\be \label{LS3-equ2.1}
\big|||P_{h+\varepsilon g}||-||P_f||\big|\leq1
\ee
holds for each $\varepsilon\in\overline{\mathbb{D}}$. In particular,
\bes
\big|||P_{h}||-||P_f||\big|\leq1.
\ees
The constant $1$  is sharp in the two estimates.
\ethm

\bpf
Suppose that $f=h+\overline{g}$ is a sense-preserving harmonic mapping in $\mathbb{D}$. Then
$h+\varepsilon g$ is a locally univalent analytic function in $\mathbb{D}$ for each $\varepsilon\in\overline{\mathbb{D}}$.
By \eqref{LS3-equ1.1}, a direct computation shows that
$$P_{h+\varepsilon g}=\frac{h''+\varepsilon g''}{h'+\varepsilon g'}
=P_h+\frac{\varepsilon \omega'}{1+\varepsilon \omega},
$$
and thus,
$$P_{h+\varepsilon g}-P_f=\frac{\varepsilon \omega'}{1+\varepsilon \omega}
+\frac{\overline{\omega}\omega'}{1-|\omega|^2}
=\frac{\varepsilon+\overline{\omega}}{1+\varepsilon\omega}\cdot\frac{\omega'}{1-|w|^2},
$$
where $\omega=g'/h'$.
Therefore, by the Schwarz-Pick lemma, we have %the inequality \eqref{LS3-equ2.1} can be deduced from
$$(1-|z|^2)\big||P_{h+\varepsilon g}(z)|-|P_f(z)|\big|\leq(1-|z|^2)\left|P_{h+\varepsilon g}(z)-P_f(z)\right|\leq
\sup_{z\in\mathbb{D}}\left|\frac{\overline{\varepsilon}
+\omega(z)}{1+\varepsilon\omega(z)}\right|\leq1
$$
for every $z\in\mathbb{D}$ and the assertion easily follows.

To show the sharpness, it suffices to consider the harmonic Koebe function  $K$ (see \cite{CS}) defined by
$$K(z)=h(z)+\overline{g(z)}=\frac{z-\frac{1}{2}z^2+\frac{1}{6}z^3}{(1-z)^3}
+\overline{\frac{\frac{1}{2}z^2+\frac{1}{6}z^3}{(1-z)^3}},\quad z\in\mathbb{D}.
$$
%which was obtained by shear construction in the horizontal direction
%of the Koebe mapping $k(z)=z/(1-z)^2$ with dilatation $\omega(z)=z$ {\rm (see \cite{CS})}.
%It is close-to-convex and always provides sharpness example in many questions
%about close-to-convex harmonic mappings.
By a direct computation, we find that
$$P_{h-g}(z)=\frac{4+2z}{1-z^2},\quad P_K(z)=\frac{5+3z}{1-z^2}-\frac{\overline{z}}{1-|z|^2},
$$
$$ P_h(z)=\frac{5+3z}{1-z^2} \quad \mbox{and} \quad P_{h+g}(z)=\frac{6+2z}{1-z^2}.
$$
It is easy to see that $||P_{h-g}||=6$ and $||P_{h}||=||P_{h+g}||=8$.
Choosing $\varepsilon=\pm1$ in \eqref{LS3-equ2.1}, it follows that $||P_{K}||=7$.
In summary, we get that
$$||P_{h-g}||+1=||P_{K}||=7=||P_{h}||-1=||P_{h+g}||-1.
$$

In addition, the sharpness can be seen from the harmonic half-plane  mapping $L$ (see \cite{CS}) defined by
$$L(z)=h(z)+\overline{g(z)}=\frac{2z-z^2}{2(1-z)^2}
+\overline{\frac{-z^2}{2(1-z)^2}},\quad z\in\mathbb{D}.
$$
%It was obtained by shear construction in the vertical direction of the function
%$l(z)=z/(1-z)$ with dilatation $\omega(z)=-z$ {\rm (see \cite{CS})}.
%It is convex and often provides sharpness example in many questions about convex harmonic mappings.
Elementary computations yield that
$$P_{h+g}(z)=\frac{2}{1-z},\quad P_L(z)=\frac{3}{1-z}-\frac{\overline{z}}{1-|z|^2},
$$
$$P_h(z)=\frac{3}{1-z}\quad \mbox{and} \quad P_{h-g}(z)=\frac{4+2z}{1-z^2}.
$$
As in the harmonic Koebe function, we obtain that
$$||P_{h+g}||+1=||P_{L}||=5=||P_{h}||-1=||P_{h-g}||-1
$$
and the proof is complete.
\epf

Obviously, the assertion \eqref{LS3-equ1.3} is true by Theorem \ref{LS3-thm2.1}.
Next we consider  more general inequalities.
%between
%$||P_{h+\varepsilon_1 \overline{g}}||$ and $|||P_{h+\varepsilon_2 g}||$,
%where $\varepsilon_1,~\varepsilon_2\in\overline{\mathbb{D}}$.

\bcor \label{LS3-cor2.1}
Let $f=h+\overline{g}$ be a sense-preserving harmonic mapping in $\mathbb{D}$.
If $||P_h||<\infty$, then for any $\varepsilon_1,~\varepsilon_2\in\overline{\mathbb{D}}$, we have the following inequalities.
\begin{enumerate}
\item The sharp inequality $\big|||P_{h+\varepsilon_1 g}||-||P_{h+\varepsilon_2 g}||\big|\leq2$ holds.

\item If $|\varepsilon_1|=|\varepsilon_2|$, then
$||P_{h+\varepsilon_1 \overline{g}}||=||P_{h+\varepsilon_2 \overline{g}}||$.
If $|\varepsilon_1|\neq|\varepsilon_2|$, then
$$\big|||P_{h+\varepsilon_1 \overline{g}}||-||P_{h+\varepsilon_2 \overline{g}}||\big|\leq|\varepsilon_1|+|\varepsilon_2|<2.
$$

\item If $|\varepsilon_1|\leq|\varepsilon_2|$, then we have the sharp inequality
$\big|||P_{h+\varepsilon_1 g}||-||P_{h+\varepsilon_2 \overline{g}}||\big|\leq1.$
If $|\varepsilon_1|>|\varepsilon_2|$, then
$$\big|||P_{h+\varepsilon_1 g}||-||P_{h+\varepsilon_2 \overline{g}}||\big|\leq1+|\varepsilon_1|+|\varepsilon_2|<3.
$$
\end{enumerate}
\ecor
\bpf
Since $||P_h||<\infty$, it follows from Theorem \ref{LS3-thm2.1} that both $||P_{h+\varepsilon g}||$
and $||P_{h+\varepsilon\overline{g}}||$ are finite for each $\varepsilon\in\overline{\mathbb{D}}$.

(1) The inequality can be easily deduced from \eqref{LS3-equ2.1} by applying the triangle inequality once.
The sharpness can be seen from the harmonic Koebe function and the harmonic half-plane  mapping.

(2) Note that $f_\varepsilon=h+\varepsilon\overline{g}$ is still a sense-preserving harmonic mapping
with dilatation $\overline{\varepsilon}\omega_f$ for any given $\varepsilon\in\overline{\mathbb{D}}$.
It follows from \eqref{LS3-equ1.1} that
$$P_{h+\varepsilon_1\overline{g}}-P_{h+\varepsilon_2 \overline{g}}
=\left( \frac{|\varepsilon_2|^2}{1-|\varepsilon_2\omega_f|^2}
-\frac{|\varepsilon_1|^2}{1-|\varepsilon_1\omega_f|^2} \right )
\overline{\omega_f}\omega'_f.
$$
Then the former part  is trivial.
The later part can be easily deduced from the Schwarz-Pick lemma and the triangle inequality.

(3) The former part is a direct consequence of \eqref{LS3-equ2.1}.
For the later part, using  (2), the former part and the triangle inequality, we have
\begin{align*}
\big|||P_{h+\varepsilon_1 g}||-||P_{h+\varepsilon_2 \overline{g}}||\big|\leq &
\big|||P_{h+\varepsilon_1 g}||-||P_{h+\varepsilon_1 \overline{g}}||\big|
+\big|||P_{h+\varepsilon_1 \overline{g}}||-||P_{h+\varepsilon_2 \overline{g}}||\big|\\
\leq & 1+|\varepsilon_1|+|\varepsilon_2|<3.
\end{align*}
The proof is complete.
\epf

Associated with Bieberbach's criterion and Yamashita's result about
convex analytic functions (see \cite[Theorem~1]{yam1999}), %(see \cite[Theorem~1]{yam1975}),
we get the following result. % as an application of \eqref{LS3-cor2.1}.

\begin{cor}\label{LS3-cor2.2}
Let $f=h+\overline{g}$ be a sense-preserving harmonic mapping in $\mathbb{D}$.
If $h+\varepsilon_1 g$ is univalent  (resp. convex) in $\mathbb{D}$ for some $\varepsilon_1\in\overline{\mathbb{D}}$,
then $||P_{h+\varepsilon\overline{g}}||<9$ (resp. $7$)
and  $||P_{h+\varepsilon g}||\leq8$ (resp. $6$) for each $\varepsilon\in\overline{\mathbb{D}}$.
Furthermore, the constants $8$ and $6$ are sharp.
Conversely, if either $||P_{h+\varepsilon_1\overline{g}}||\geq9$
(resp. $7$) or $||P_{h+\varepsilon_2 g}||>8$ (resp. $6$)
for some $\varepsilon_1,~\varepsilon_2\in\overline{\mathbb{D}}$,
then $h+\varepsilon g$ is not univalent (resp. convex) in $\mathbb{D}$ for any $\varepsilon\in\overline{\mathbb{D}}$.
\end{cor}

The harmonic Koebe function $K=h_K+\overline{g_K}$ and the harmonic half-plane mapping $L=h_L+\overline{g_L}$ still show its sharpness
in the corresponding cases because $h_K(z)-g_K(z)=\frac{z}{(1-z)^2}$ is univalent in $\ID$ and $h_L(z)+g_L(z)=\frac{z}{1-z}$ is univalent and
convex in $\mathbb{D}$, respectively.

\section{Pre-Schwarzian norm and uniformly hyperbolic radius} \label{LS3-sec3}

It is natural to ask whether there exists a generalization of Bieberbach's criterion
for univalent harmonic mappings.
%Since the PSD preserves affine invariance, it suffices to  consider only functions in $\mathcal{H}$ or $\mathcal{H}_0$.
Let
$$\mathcal{S}_H=\left \{f\in\mathcal{H}:\,f(z)=h(z)+\overline{g(z)}
=\sum_{n=1}^{\infty}a_nz^n+\sum_{n=1}^{\infty}\overline{b}_n\overline{z}^n
~\mbox{is univalent in} ~\mathbb{D}\right \}
$$
and
$\mathcal{S}_H^0=\mathcal{S}_{H}\cap\mathcal{H}_0.$
Set
$$\alpha=\sup_{f\in \mathcal{S}_H}|a_2|
\quad \mbox{and} \quad \alpha_0=\sup_{f\in \mathcal{S}_H^0}|a_2|.
$$
Clunie and Sheil-Small \cite{CS} showed that if $f=h+\overline{g}\in\mathcal{S}_H$,
then $||P_{h}||\leq2(\alpha+1)$,
$\alpha_0<12172$ and $\alpha_0\leq\alpha\leq\alpha_0+1/2$.
They conjectured that $\alpha_0\leq5/2$, which has a special significance
in many extremal problems for harmonic mappings.
The estimate of $\alpha_0$ was improved (see \cite[p.~96]{dur2004} and \cite[Theorem~10]{she}).
Now the best known upper bound for  $\alpha_0$ is in \cite{AAP}.

However, for certain geometric subfamilies of $\mathcal{S}_H$,
 %the family of univalent harmonic mappings \textcolor[rgb]{1.00,0.00,0.00}{in $\mathcal{H}$,}
we have some precise coefficient estimates.
For example, for the families $\mathcal{K}_{H}$ and $\mathcal{C}_{H}$ of
convex and close-to-convex harmonic mappings in $\mathbb{D}$, respectively.
We note that $\mathcal{K}_{H}\subseteq \mathcal{C}_{H} \subseteq \mathcal{S}_{H}$.
Set $\mathcal{K}_H^0=\mathcal{K}_{H}\cap \mathcal{H}_0$
and $\mathcal{C}_H^0=\mathcal{C}_{H}\cap \mathcal{H}_0$.
For these special families,  we know (see \cite{CS} and \cite{WLZ}):
$$\sup_{f\in \mathcal{K}_H^0}|a_2|=\frac{3}{2},\quad
\sup_{f\in \mathcal{K}_H}|a_2|=2,
\quad \sup_{f\in \mathcal{C}_H^0}|a_2|=\frac{5}{2}\quad \mbox{and}
\quad \sup_{f\in \mathcal{C}_H}|a_2|=3.
$$
Therefore, the sharp estimate $||P_f||\leq5$ is obtained for all $f\in\mathcal{K}_{H}$ (see \cite[Theorem~4]{HM2015}).
On the other hand, based on further research on affine and linear invariant
locally univalent harmonic mappings, Graf in \cite[Theorem~1]{gra} obtained that $||P_{f}||\leq7$
for $f\in\mathcal{C}_{H}$ and $||P_{f}||\leq2(\alpha_0+1)$ for $f\in\mathcal{S}_{H}$.

In this section, we will first re-certify the above partial results  concerning the PSN
as a direct consequence of our present study on ULU harmonic mappings. For the convenience of
the reader, we include the proof here since it follows by a direct computation. Note that
the PSN is in general not linear invariant.

\bthm\label{LS3-thm3.1}
Let $f=h+\overline{g}$ be a sense-preserving and {\rm ULU} harmonic mapping in $\mathbb{D}$.
Then we have
\begin{equation} \label{LS3-equ3.1}
(1-|z|^2)|P_{h}(z)|\leq2(\alpha/t+|z|)\quad \mbox{and}\quad (1-|z|^2)|P_{f}(z)|\leq2(\alpha_0/t+|z|)
\end{equation}
for  every $z\in\mathbb{D}$, where
$$t= \left\{\begin{array}{cl}
\ds \frac{e^{\rho(f)}-1}{e^{\rho(f)}+1} & \mbox{ if $\rho(f)<\infty$},\\
1& \mbox{ if $\rho(f)=\infty$.} \end{array}\right .
$$
%and $t$ equals 1 if $\rho(f)=\infty$.
 In particular, if $f$ is univalent in $\mathbb{D}$, then
\begin{equation*} %\label{LS3-equ4.1}
||P_{h}||\leq2(\alpha+1)\quad \mbox{and}\quad ||P_{f}||\leq2(\alpha_0+1).
\end{equation*}
\ethm
\bpf
Suppose that $f=h+\overline{g} \in {\rm ULU}$. Then
$f$ is univalent in each hyperbolic disk $d_h(z,\rho(f))$ for every $z\in\mathbb{D}$.
Fix $z\in\mathbb{D}$ and let $\phi(\zeta)=\frac{t\zeta+z}{1+t\overline{z}\zeta}$ ($\zeta\in\mathbb{D}$),
where $t$ is defined as above.
Using the Koebe transformation, we get that
\begin{align*}
F_1(\zeta)=&\frac{(f\circ\phi)(\zeta)-(f\circ\phi)(0)}{(f\circ\phi)_\zeta(0)}\\
=&\frac{h(\phi(\zeta))-h(z)}{th'(z)(1-|z|^2)}+\overline{\frac{g(\phi(\zeta))-g(z)}{t\overline{h'(z)}(1-|z|^2)}}\\
=& H_1(\zeta)+\overline{G_1(\zeta)}
\end{align*}
and $F_1\in \mathcal{S}_{H}$. A simple computation yields that
$$|H''_1(0)|=t\left|(1-|z|^2)\frac{h''(z)}{h'(z)}-2\overline{z}\right|\leq 2\alpha,
$$
which implies the first inequality in \eqref{LS3-equ3.1}. Using the affine change, we have that
\begin{align*}
F_2(\zeta)=&\frac{F_1(\zeta)-\overline{b_1F_1(\zeta)}}{1-|b_1|^2}\\
=&\frac{H_1(\zeta)-\overline{b_1}G_1(\zeta)}{1-|b_1|^2}+\frac{\overline{G_1(\zeta)-b_1H_1(\zeta)}}{1-|b_1|^2}\\
=& H_2(\zeta)+\overline{G_2(\zeta)}
\end{align*}
%where $F_2\in\mathcal{S}_{H}^0$
%and $b_1=G'_1(0)=g'(z)/\overline{h'(z)}$.
and $F_2\in\mathcal{S}_{H}^0$, where $b_1=G'_1(0)=g'(z)/\overline{h'(z)}$.
Again, a straightforward computation shows that
\begin{align*}
|H''_2(0)|=&\frac{|H''_1(0)-\overline{b_1}G''_1(0)|}{1-|b_1|^2}\\
=&t\left|(1-|z|^2)\frac{h''(z)\overline{h'(z)}-g''(z)\overline{g'(z)}}{|h'(z)|^2-|g'(z)|^2}-2\overline{z}\right|\\
=& t|(1-|z|^2)P_f(z)-2\overline{z}|\leq2\alpha_0,
\end{align*}
which implies the second inequality in \eqref{LS3-equ3.1}.
\epf

%\begin{cor}\label{LS3-cor3.1}
%Let $f=h+\overline{g}$ be a sense-preserving harmonic mapping in $\mathbb{D}$.
%If $f$ is convex (resp. close-to-convex) in $\mathbb{D}$, then we have
%$$||P_{f}||\leq5~~~(resp.~7) \quad \mbox{and}\quad ||P_{h+\varepsilon g}||\leq6~~~(resp.~8)
%$$
%for each $\varepsilon\in\overline{\mathbb{D}}$.
%All estimates are sharp.
%\end{cor}
%
%The proof of the first parter in Corollary \ref{LS3-cor3.1} is similar to Theorem \ref{LS3-thm3.1}
%and the last parter follows from \eqref{LS3-equ2.1}.
%Moreover, the harmonic Koebe function and harmonic half-plane  mapping show the sharpness
%for close-to-convex and convex harmonic mappings, respectively.

Next we consider stable harmonic univalent (resp. convex) mappings.
A sense-preserving harmonic mapping $f=h+\overline{g}$  in $\mathbb{D}$ is called SHU
(resp.  SHC)  if $h+\lambda\overline{g}$
is univalent (resp. convex)  in $\mathbb{D}$  for every $|\lambda|=1$.
The following result has some similarities with the classical estimate of the SD
for SHU and SHC mappings in \cite[Theorem~2]{CHM}, but the method of proof is
different and so can also be adapted to prove \cite[Theorem~2]{CHM}.

\begin{thm}\label{LS3-thm3.2}
Let $f=h+\overline{g}$ be a sense-preserving harmonic mapping in $\mathbb{D}$.
If $f$ is {\rm SHU} (resp. {\rm SHC}), then we have
$$||P_{h+\varepsilon g}||\leq6~~~\mbox{(resp.~{\rm 4})}
\quad \mbox{and} \quad ||P_{h+\varepsilon \overline{g}}||\leq6~~~\mbox{(resp.~{\rm 4})}
$$
for each $\varepsilon\in\overline{\mathbb{D}}$.
All estimates are sharp.
\end{thm}
\bpf
If $f=h+\overline{g}$ is SHU (resp. SHC) in $\mathbb{D}$,
then both $h+\varepsilon g$ and $h+\varepsilon \overline{g}$ are univalent (resp. convex) in $\mathbb{D}$
for each $\varepsilon\in\overline{\mathbb{D}}$ (see \cite{HM2013}).
It follows from Bieberbach's criterion (resp. \cite[Theorem~1]{yam1999})
that $||P_{h+\varepsilon g}||\leq6$ (resp. $4$) for each $\varepsilon\in\overline{\mathbb{D}}$.

Fix $\varepsilon\in\overline{\mathbb{D}}$ and let $f_\varepsilon=h+\varepsilon\overline{g}$.
For all $z_0\in\mathbb{D}$, it follows from \cite[Lemma~1]{HM2015} that
$P_{f_\varepsilon}(z_0)=P_{h-\varepsilon\overline{\omega(z_0)}g}(z_0)$ and thus,
$$(1-|z_0|^2)|P_{f_\varepsilon}(z_0)|=(1-|z_0|^2)|P_{h-\varepsilon\overline{\omega(z_0)}g}(z_0)|,
$$
where $\omega=g'/h'$.
This implies that $||P_{f_\varepsilon}||\leq\sup_{\lambda\in\overline{\mathbb{D}}}||P_{h+\lambda g}||$
and the assertion follows.

%in particular, $h$ itself is convex (univalent).
%Thus, for any given $\varepsilon\in\overline{\mathbb{D}}$, it is clearly that
%$h+\lambda(\overline{\varepsilon}g)$ is convex (resp. univalent) for every $|\lambda|=1$.
%It follows from \cite[Theorem~3.1]{HM2013} (resp. \cite[Proposition~2.1]{HM2013})
%that $h+\varepsilon\overline{g}$ is also SHC (resp. SHU).
%Note that PSD preserves affine invariance property and convexity (resp. univalence)
%under the affine harmonic mapping.
%Without loss of generality, we can restricted the SHC (resp. SHU) in $\mathcal{S}_{H}^0$.
%It follows from \cite[Theorem~8.1]{HM2013} and \cite[Proposition~8.2]{HM2013}   that
%$$\sup_{f\in \mathcal{S}^0_{H}\cap \text{SHC}}|a_2|=1 \quad \text{and} \quad
%\sup_{f\in \mathcal{S}^0_{H}\cap \text{SHU}}|a_2|=2,
%$$
%respectively. By the similar proof of Theorem \ref{LS3-thm3.1}, we can get
%$||P_{h+\varepsilon \overline{g}}||\leq4$ (resp. $6$) for each $\varepsilon\in\overline{\mathbb{D}}$.

To show that all estimates are sharp, it is enough to consider the analytic functions
$$k(z)=\frac{z}{(1-z)^2}\quad\mbox{and}\quad l(z)=\frac{1+z}{1-z}
$$
that belong to the families of SHU and SHC mappings with $||P_k||=6$ and $||P_l||=4$, respectively.
\epf

Combining Corollary \ref{LS3-cor2.1} and Theorem \ref{LS3-thm3.1} (resp. Theorem \ref{LS3-thm3.2}),
we can obtain a few similar results as that of Corollary \ref{LS3-cor2.2} for univalent harmonic mappings (resp. SHU and SHC mappings).
However, we do not include these statements here.  Below we consider the converse of Theorem \ref{LS3-thm3.1}.
%As a preparation, by Lemma 1 and Theorem 2 in \cite{son}, we have the following result.
%
%\blem \label{LS3-lem3.1}
%Let $h$ be a locally univalent analytic function in $\mathbb{D}$.
%If $||P_{h}||\leq M$, then $h$ is univalent in $D_h(z,2\tanh^{-1}1/(8t))$
%for each $z\in\mathbb{D}$, where $t=\max\{M/2,1/4\}$.
%\elem

\bthm\label{LS3-thm3.3}
Let $f=h+\overline{g}$ be a sense-preserving harmonic mapping in $\mathbb{D}$.
If $||P_{f}||\leq M$, then $f$ is univalent in the hyperbolic disk $D_h(z,t)$ for each $z\in\mathbb{D}$.
Consequently, $f$ is {\rm ULU} in $\mathbb{D}$ and its uniformly hyperbolic radius $\rho(f)$ is no less than $t$.
Here $t=2\tanh^{-1}\left(1/(8(M+1))\right)$.
\ethm
\bpf
Fix $\varepsilon\in \mathbb{D}$. Let $f_\varepsilon=f+\varepsilon\overline{f}=\phi_\varepsilon+\overline{\psi_\varepsilon}$, where
$\phi_\varepsilon=h+\varepsilon g$ and $\psi_\varepsilon=g+\overline{\varepsilon}h$.
%Since the dilatation of $f_\varepsilon$ equals to $\phi\circ\omega_f$  for each $\varepsilon\in\mathbb{D}$,
%where $\phi(z)=\frac{z+\varepsilon}{1+\overline{\varepsilon}z}\in{\rm Aut}(\mathbb{D})$,
%It follows from \eqref{LS3-equ1.5} that $f_\varepsilon$ is still
%a sense-preserving harmonic mapping for each $\varepsilon\in\mathbb{D}$.
By the hypothesis,  \eqref{LS3-equ1.2} and \eqref{LS3-equ2.1}, we get that
$$||P_{\phi_\varepsilon}||\leq||P_{f_\varepsilon}||+1=||P_{f}||+1\leq M+1.
$$
It follows from  \cite[Theorem~2]{son} that $\phi_\varepsilon$ is univalent in $D_h(z,t)$ for each $z\in\mathbb{D}$,
where $t$ as the above.
%$$t=2\tanh^{-1}\left(1/(8(M+1))\right).
%$$
By Hurwtiz's theorem, we know that for each $z\in\mathbb{D}$,
$h+\lambda g$ is univalent in $D_h(z,t)$ for every $|\lambda|=1$.
Therefore, it follows from \cite[ Corollary~2.2]{HM2013} that $f$ is univalent in $D_h(z,t)$ for each $z\in\mathbb{D}$.
This ends the proof.
\epf

\section{Stable geometric properties of ULU analytic and harmonic mappings} \label{LS3-sec4}
In this section, we will show a great number of equivalent conditions for sense-preserving and ULU harmonic mappings in $\mathbb{D}$.
First we will introduce some notations.
Let $f=h+\overline{g}$ be a sense-preserving harmonic mapping in $\mathbb{D}$. Set
$$\rho^*(f)=\inf_{z\in\mathbb{D}}\left\{\sup_{\rho_z>0}\{\rho_z:\, f ~\mbox{is convex in}~ D_h(z,\rho_z)\}\right\}.
$$
If $\rho^*(f)>0$, then we say that $f\in  {\rm ULC}$.
%Clearly, an ULC analytic function is a special  ULC harmonic mapping.
%The {\it Schwarzian derivative and  norm} of $f$ were introduced in \cite{HM2015} (see also \cite{CDO}) and defined by
The SD and the SN of $f$ were investigated in details by Hern\'{a}ndez and Mart\'{\i}n  \cite{HM2015} (see also \cite{CDO}) and
they were defined by
$$S_f=S_h+\frac{\overline{\omega}}{1-|\omega|^2}\left(\frac{h''}{h'}\omega'
-\omega''\right)-\frac{3}{2}\left(\frac{\overline{\omega}\omega'}{1-|\omega|^2}\right)^2%, \quad \mbox{for}\quad z\in\mathbb{D},
$$
and
$$||S_{f}||=\sup_{z\in\mathbb{D}}(1-|z|^2)^2|S_{f}(z)|,
$$
respectively, where $S_h$ is the classical Schwarzian derivative of a
 locally univalent function $h$ defined by %$P_f$ is defined as \eqref{LS3-equ1.1} and
$$S_h=%(P_h)'-\frac{1}{2}(P_h)^2=
\frac{h'''}{h'}-\frac{3}{2}\left(\frac{h''}{h'}\right)^2.%,
%\quad \mbox{for} \quad z\in\mathbb{D},
$$
%which is $S_h$ is the classical Schwarzian derivative of locally univalent analytic function $h$ in $\mathbb{D}$.
If $g$ is a constant, then it is clear that $S_f=S_h$ and $||S_f||=||S_h||$.
%The  Schwarzian norm of $h$ is defined by
%$$||S_{h}||=\sup_{z\in\mathbb{D}}|S_{h}(z)|(1-|z|^2)^2,
%$$
%where $S_h$ is its Schwarzian derivative, namely,
%$$S_h=(P_h)'-\frac{1}{2}(P_h)^2=
%\frac{h'''}{h'}-\frac{3}{2}\left(\frac{h''}{h'}\right)^2,
%\quad z\in\mathbb{D}.
%$$
Analogous to some features of  the PSN, if $f$ is a univalent analytic function in $\mathbb{D}$,
we have the sharp inequality $||S_{f}||\leq6$.
Conversely, if  $||S_{f}||\leq2$ for a locally univalent analytic function $f$
in $\mathbb{D}$, then, according to Krauss-Nehari's criterion, $f$ is univalent in $\mathbb{D}$
and the constant $2$ is sharp (see \cite{kra,neh}).
%Furthermore, these constants are sharp.
There are some criteria for the univalence of harmonic mappings in terms of the SN (see \cite{gra,HM2015,HM2015-1}),
but these results are not sharp.
%However, even if for convex and close-to-convex harmonic mapping $f$ in $\mathbb{D}$,
%the best known upper bound for $||S_f||$ is $6$ and $26.368$, respectively (see \cite[Theorem~3]{gra} and \cite[Proposition~2]{HM2015}).
%Conversely, if  $f$ is a locally univalent harmonic mapping in $\mathbb{D}$ with $||S_{f}||\leq \delta$  for sufficiently $\delta>0$,
% then $f$ is univalent in $\mathbb{D}$ (see \cite[Theorem~2.2]{HM2015-1}).

Next, we present equivalent conditions for sense-preserving and ULU harmonic mappings
in $\mathbb{D}$ based on %the results in \cite{HM2015} and
the following result.

\blem   {\rm (\cite[p.~44]{dur1983} and \cite[Theorem~2]{yam1977})} \label{LS3-lem4.1}
Let $f$ be a locally univalent analytic function  in $\mathbb{D}$.
Then the following are equivalent.
\begin{enumerate}
\item $f\in {\rm ULU}$; %$f$ is  {\rm ULU} in $\mathbb{D}$;
\item $f\in  {\rm ULC}$; %$f$ is  {\rm ULC} in $\mathbb{D}$;
\item $||P_f||<\infty$;
\item $||S_f||<\infty$;
\item There exists a constant $m>0$, and a univalent analytic function $F$  in $\mathbb{D}$
such that $f'=(F')^m$.
\end{enumerate}
\elem

%Chuaqui et al. \cite{CDO} seem to the first
%who extend the definition mentioned above to the harmonic mapping $f$ defined in a simply connected domain
%with dilatation $\omega_f=q^2$ for some analytic function $q$.
%Obviously, such definition is required to fulfill extra condition
%so that it can't be applied to any locally univalent harmonic mapping.

To describe our results, we introduce the following abbreviations analogous to the paper \cite{HM2013}.
Below, let $f=h+\overline{g}$ be a sense-preserving harmonic mapping in $\mathbb{D}$
and $\varepsilon_1,~\varepsilon_2\in\overline{\mathbb{D}}$ with $\varepsilon_1\neq\varepsilon_2$.
If $h+\varepsilon_1\overline{g}$ is ULU (resp. ULC) in $\mathbb{D}$
if and only if $h+\varepsilon_2\overline{g}$ is ULU (resp. ULC) in $\mathbb{D}$,
then we say that $f$ is SHULU (resp.  SHULC).
Similarly, if $h+\varepsilon_1g$ is ULU (resp. ULC) in $\mathbb{D}$
if and only if $h+\varepsilon_2g$ is ULU (resp. ULC) in $\mathbb{D}$,
then we say that $h+g$ is SAULU
(resp. SAULC).
%Next we introduce the definition of stable bounded  for pre-Schwarzian and Schwarzian norms.
If $||P_{h+\varepsilon_1\overline{g}}||$ (resp. $||S_{h+\varepsilon_1\overline{g}}||$) is bounded
if and only if $||P_{h+\varepsilon_2\overline{g}}||$ (resp. $||S_{h+\varepsilon_2\overline{g}}||$) is bounded,
then we say that $f$ has SBHPSN
(resp. SBHSN).
Similarly, if $||P_{h+\varepsilon_1g}||$ (resp. $||S_{h+\varepsilon_1g}||$) is bounded
if and only if $||P_{h+\varepsilon_2g}||$ (resp. $||S_{h+\varepsilon_2g}||$) is bounded,
then we say that $h+g$ has SBAPSN
(resp. SBASN).

\bthm \label{LS3-thm4.1} {\rm (Equivalent conditions)}
Let $f=h+\overline{g}$ be a sense-preserving harmonic mapping in $\mathbb{D}$.
Then the following conditions are equivalent.

\begin{enumerate}
\item $h+g$ is {\rm SAULU};

\item $h+g$ is {\rm SAULC};

\item $h+g$ has {\rm SBAPSN};

\item $h+g$ has  {\rm SBASN};

\item  For any two points $\varepsilon_1,~\varepsilon_2\in\overline{\mathbb{D}}$ with $\varepsilon_1\neq\varepsilon_2$,
 there exists a constant $m_1>0$, and a univalent analytic function $F_1$ such that $(h+\varepsilon_1g)'=(F'_1)^{m_1}$
if and only if there exists a constant $m_2>0$, and a univalent analytic function $F_2$ such that
$$(h+\varepsilon_2g)'=(F'_2)^{m_2};
$$

\item $f$ is {\rm SHULU};

\item $f$ is {\rm SHULC};

\item $f$ has {\rm SBHPSN};

\item $f$ has {\rm SBHSN}.
\end{enumerate}
\ethm
\bpf
%The scheme of the proof is to show that %$(3)\Leftrightarrow(8)$.
%%Then we prove that
%$(1)\Leftrightarrow(2)\Leftrightarrow(3)\Leftrightarrow(4)\Leftrightarrow(5)\Leftrightarrow(8)$
%and  $(6)\Leftrightarrow(7)\Leftrightarrow(8)\Leftrightarrow(9)$.
To simplify the proof, we use the equivalent diagram below.
%Thus, we have that $(3)\Leftrightarrow(8)$.
If we apply  Lemma \ref{LS3-lem4.1} to $h+\varepsilon_1g$ and $h+\varepsilon_2g$,
we see that (Ai) and (Bi)~$(i=1,2,3,4)$ hold.
On the other hand, (A5), (AB) and (B5) are the direct consequences of Theorem \ref{LS3-thm2.1} and Corollary \ref{LS3-cor2.1}.
%Then by the bridge of (A5), (AB) and (B5),
Clearly, the following implications are easy to obtain
$$(1)\Leftrightarrow(2)\Leftrightarrow(3)\Leftrightarrow(4)\Leftrightarrow(5)\Leftrightarrow(8).
$$

\noindent
{\small
$\renewcommand{\arraystretch}{1.6}
\begin{array}[c]{ccccccc}
||S_{h+\varepsilon_1g}||<\infty&\stackrel{A1}{\Longleftrightarrow}&(h+\varepsilon_1g)'
=(F'_1)^{m_1}&&(h+\varepsilon_2g)'=(F'_2)^{m_2}
& \stackrel{B1}{\Longleftrightarrow}\|S_{h+\varepsilon_2g}||<\infty \\  [-0.5ex]
\Updownarrow\scriptstyle{A2}&&&&&\Updownarrow\scriptstyle{B2}& \\ [-0.5ex]
h+\varepsilon_1g\in\mbox{ULU}&\stackrel{A3}{\Longleftrightarrow}&h+\varepsilon_1g
\in\mbox{ULC}&&h+\varepsilon_2g\in\mbox{ULC}
&\stackrel{B3}{\Longleftrightarrow}  h+\varepsilon_2g\in\mbox{ULU}\\ [-0.5ex]
\Updownarrow\scriptstyle{A4}&&&&& \Updownarrow\scriptstyle{B4}&\\ [-0.5ex]
||P_{h+\varepsilon_1g}||<\infty&\stackrel{A5}{\Longleftrightarrow}&||P_{h+\varepsilon_1\overline{g}}||<\infty
&\stackrel{AB}\Longleftrightarrow& ||P_{h+\varepsilon_2\overline{g}}||<\infty&
\stackrel{B5}{\Longleftrightarrow}||P_{h+\varepsilon_2g}||<\infty \\ [-0.5ex]
&&\Updownarrow\scriptstyle{A6}&&\Updownarrow\scriptstyle{B6}&&\\ [-0.5ex]
||S_{h+\varepsilon_1\overline{g}}||<\infty&\stackrel{A7}{\Longleftrightarrow}
& h+\varepsilon_1\overline{g}\in\mbox{ULU}&&h+\varepsilon_2\overline{g}
\in\mbox{ULU}& \stackrel{B7}{\Longleftrightarrow}||S_{h+\varepsilon_2\overline{g}}||<\infty \\ [-0.5ex]
&&\Updownarrow\scriptstyle{A8}&&\Updownarrow\scriptstyle{B8}&&\\ [-0.5ex]
&&h+\varepsilon_1\overline{g}\in\mbox{ULC}&&h+\varepsilon_2\overline{g}\in\mbox{ULC}&&
\end{array}
$
}

To complete the proof, we need to show that $(6)\Leftrightarrow(7)\Leftrightarrow(8)\Leftrightarrow(9)$.
If we apply Theorems \ref{LS3-thm3.1} and \ref{LS3-thm3.3} to $h+\varepsilon_1\overline{g}$ and $h+\varepsilon_2\overline{g}$,
then we obtain the inclusions (A6) and (B6).

From \cite[Theorem~7]{HM2015},  (A7) and (B7) follow.

To prove (A8) and (B8), it suffices to show that each $f=h+\overline{g}\in {\rm ULU}$ also belongs to ULC.
To do this, let us assume that $f=h+\overline{g}$ is ULU in $\mathbb{D}$. Then  $M=||P_f||<\infty$ and thus,
$\sup_{\varepsilon\in\overline{\mathbb{D}}}||P_{h+\varepsilon g}||\leq M+1$ by \eqref{LS3-equ2.1}.
Following the proof and notations of Theorem \ref{LS3-thm3.3}, we see that for each $z\in\mathbb{D}$,
$h+\lambda g$ is convex in $D_h(z,(2-\sqrt{3})t)$ for every $|\lambda|=1$ by the classical result on the radius of convexity (see \cite[p.~44]{dur1983}),
where $t=2\tanh^{-1}(1/(8(M+2)))$.
It follows from \cite[Theorem~3.1]{HM2013} that $f$ is convex in the hyperbolic disk
$D_h(z,(2-\sqrt{3})t)$ for each $z\in\mathbb{D}$, which means that $f$ is ULC in $\mathbb{D}$.

Again, by the bridge (AB), we prove that  $(6)\Leftrightarrow(7)\Leftrightarrow(8)\Leftrightarrow(9)$. This completes the proof.
%
%\noindent
%{\small
%$\renewcommand{\arraystretch}{1.6}
%\begin{array}[c]{ccccccc}
%||S_{h+\varepsilon_1g}||<\infty&\stackrel{A1}{\Longleftrightarrow}&(h+\varepsilon_1g)'
%=(F'_1)^{m_1}&&(h+\varepsilon_2g)'=(F'_2)^{m_2}
%& \stackrel{B1}{\Longleftrightarrow}\|S_{h+\varepsilon_2g}||<\infty \\  [-0.5ex]
%\Updownarrow\scriptstyle{A2}&&&&&\Updownarrow\scriptstyle{B2}& \\ [-0.5ex]
%h+\varepsilon_1g\in\mbox{ULU}&\stackrel{A3}{\Longleftrightarrow}&h+\varepsilon_1g
%\in\mbox{ULC}&&h+\varepsilon_2g\in\mbox{ULC}
%&\stackrel{B3}{\Longleftrightarrow}  h+\varepsilon_2g\in\mbox{ULU}\\ [-0.5ex]
%\Updownarrow\scriptstyle{A4}&&&&& \Updownarrow\scriptstyle{B4}&\\ [-0.5ex]
%||P_{h+\varepsilon_1g}||<\infty&\stackrel{A5}{\Longleftrightarrow}&||P_{h+\varepsilon_1\overline{g}}||<\infty
%&\stackrel{AB}\Longleftrightarrow& ||P_{h+\varepsilon_2\overline{g}}||<\infty&
%\stackrel{B5}{\Longleftrightarrow}||P_{h+\varepsilon_2g}||<\infty \\ [-0.5ex]
%&&\Updownarrow\scriptstyle{A6}&&\Updownarrow\scriptstyle{B6}&&\\ [-0.5ex]
%||S_{h+\varepsilon_1\overline{g}}||<\infty&\stackrel{A7}{\Longleftrightarrow}
%& h+\varepsilon_1\overline{g}\in\mbox{ULU}&&h+\varepsilon_2\overline{g}
%\in\mbox{ULU}& \stackrel{B7}{\Longleftrightarrow}||S_{h+\varepsilon_2\overline{g}}||<\infty \\ [-0.5ex]
%&&\Updownarrow\scriptstyle{A8}&&\Updownarrow\scriptstyle{B8}&&\\ [-0.5ex]
%&&h+\varepsilon_1\overline{g}\in\mbox{ULC}&&h+\varepsilon_2\overline{g}\in\mbox{ULC}&&
%\end{array}
%$
%}
\epf

\begin{rems}  In the remarks below, let $f=h+\overline{g}$ be sense-preserving  in $\mathbb{D}$.

\begin{enumerate}
\item[(1)] The pre-Schwarzian norm $||P_{g}||$ and the Schwarzian norm $||S_{g}||$ can be unbounded
even if  $f$ and $g$ are univalent and locally univalent in $\mathbb{D}$, respectively.
For example, let
$$f_n(z)=h_n(z)+\overline{g_n(z)}=z-1+\overline{\lambda(z-1)^n}~~(n\geq2 ~\mbox{and}~ 0<|\lambda|<1/(n2^{n-1})).
$$
It is easy to see that $f_n$ is sense-preserving and univalent in $\mathbb{D}$
 and $g_n$ is locally univalent in $\mathbb{D}$ for any $n\geq2$.
However, we have that
$$||P_{g_n}||=\sup_{z\in\mathbb{D}}(1-|z|^2)\frac{n-1}{|1-z|}=2(n-1)\rightarrow\infty
$$
and
$$||S_{g_n}||=\sup_{z\in\mathbb{D}}(1-|z|^2)^2\frac{n^2-1}{2|1-z|^2}=2(n^2-1)\rightarrow\infty
$$
as $n\rightarrow\infty$.

\item[(2)] On one hand, the dilatation of $f$ can be expressed as square of certain analytic function if both $h$ and $g$ are  ULU in $\mathbb{D}$.
It follows from Lemma \ref{LS3-lem4.1} that $||P_{h}||$ and $||P_{g}||$  are bounded.
Let
$k=\max\{||P_{h}||,||P_{g}|| \}+1$ and set
$$h_1(z)=\int_0^z(h'(\zeta))^\frac{1}{2k}d\zeta\quad \mbox{and}\quad g_1(z)=\int_0^z(g'(\zeta))^\frac{1}{2k}d\zeta
$$
in the proof of  \cite[Theorem~2]{yam1977}.
Note that $h_1$ and $g_1$ are analytic and univalent in $\mathbb{D}$ such that
$h'=(h'_1)^{2k}$ and $g'=(g'_1)^{2k}$.
Thus, we have that $\omega_f=g'/h'=(g'_1/h'_1)^{2k}$.
Furthermore, if $f$ is univalent in $\mathbb{D}$, then $f$ can be lifted to a regular minimal surface
given by conformal (or isothermal) parameters in $\mathbb{D}$.

\item[(3)] On the other hand, the function $f$, with the dilatation $\omega_f=q^2$,
for some analytic function $q$ may not belong to ULU. For instance, let
$$f(z)=h(z)+\overline{g(z)}
=e^{-2\frac{z+1}{z-1}}+\overline{\frac{z-5}{z-1}}, \quad z\in\mathbb{D}.
$$
A simple computation infers that
$$\omega_{f}(z)=\frac{g'(z)}{h'(z)} =\left (e^{\frac{z+1}{z-1}}\right )^{2}
$$
so that $|\omega_{f}(z)|<1$ in $\mathbb{D}$ and thus,  $f$ is sense-preserving in $\mathbb{D}$.
However,
$$||P_{h}||=\sup_{z\in\mathbb{D}}(1-|z|^2)\frac{|6-2z|}{|1-z|^2}=\infty,
$$
which implies that $f$ is not ULU in $\mathbb{D}$  by Theorem \ref{LS3-thm4.1}.

\item[(4)] If the analytic part $h$ is univalent in $\mathbb{D}$,
then $f$  is certainly ULU in $\mathbb{D}$ by  Corollary \ref{LS3-cor2.2} and Theorem \ref{LS3-thm4.1}.
However, the above example shows that even if the co-analytic part $g$ is univalent in $\mathbb{D}$, $f$ may not belong to ULU.
\end{enumerate}
\end{rems}

\section{Some precise examples}\label{LS3-sec5}

In this section, we consider a family of harmonic mappings and compute their PSNs and then discuss the univalency of
the corresponding mapping. We next introduce
\begin{align} \label{LS3-equ5.1}
F_{a,b,\theta}(z)&=H_{a,b}(z)+\overline{G_{a,b,\theta}(z)}, \quad
e^{-i\theta}G_{a,b,\theta}(z):=   G_{a,b}(z) =  H_{a+1,b}(z)-H_{a,b}(z) ,
\end{align}
%\begin{align}
%F_{a,b,\theta}(z)&=H_{a,b}(z)+\overline{G_{a,b,\theta}(z)} \label{LS3-equ5.1}\\%&=&\int_0^z\frac{(1+t)^a}{(1-t)^b}dt
%%+\overline{e^{i\theta}\int_0^zt\frac{(1+t)^a}{(1-t)^b}dt}\\
%\nonumber &=H_{a,b}(z)+\overline{e^{i\theta}G_{a,b}(z)}\\
%\nonumber &=H_{a,b}(z)+\overline{e^{i\theta}(H_{a+1,b}(z)-H_{a,b}(z))},
%\end{align}
where $a,~b,~\theta\in\mathbb{R}$ and
\begin{equation}\label{LS3-equ5.2}
H_{a,b}(z)=\int_0^z\frac{(1+t)^a}{(1-t)^b}dt.
\end{equation}
If $a=b$, we denote $H_{a,a}$ by $H_a$.
Clearly, $H_{a,b}\in \mathcal{A}$ and $H_{a,b}(z)=-H_{-b,-a}(-z)$.
Therefore, it is easy to see that $F_{a,b,\theta}\in \mathcal{H}_0$ with dilatation $\omega(z)=e^{i\theta}z$ and
\begin{equation}\label{LS3-equ5.3}
F_{a,b,\theta}(z)=-F_{-b,-a,\theta+\pi}(-z),\quad z\in\mathbb{D}.
\end{equation}
In general, computing the PSN and verifying the univalence of a given harmonic mapping are not so easy.
Below, we also try to give partial answers to this issue.
Moreover, as a byproduct of our investigation, we present some sharp inequalities in Section \ref{LS3-sec6}
and give certain properties of the family $\mathcal{B}_H(\lambda)$ ($\lambda\geq1$). In the following results, we use the
following well-known facts: If $h$ is a normalized (i.e. $h(0)=h'(0)-1=0$) analytic function in $\ID$ satisfying the condition
$${\rm Re}\left(1+\frac{zh''(z)}{h'(z)}\right)>-\frac{1}{2}
$$
for $|z|<1$, then $h$ is convex in some direction and hence it is close-to-convex (univalent) in the unit disk.
For details and its importance see \cite{Hiroshi-Samy-2014}.

\begin{prop}\label{LS3-prop5.1}
For the functions $H_{a,b}$ and $H_a$ defined by \eqref{LS3-equ5.2}, we have the following properties:
\begin{enumerate}
\item $\|P_{H_{a,b}}\|=2\max\{|a|,|b|\}$.
Thus, if $\max\{|a|,|b|\}\leq 1/2$, then the functions $H_{a,b}$ are univalent in $\mathbb{D}$.
If $\max\{|a|,|b|\}>3$, then the functions $H_{a,b}$ are not univalent in $\mathbb{D}$.
 %by Bieberbach's criterion.

\item If $\min\{|a|,|b|\}+|a-b|\leq1$, then  the functions $H_{a,b}$ are close-to-convex and univalent in $\mathbb{D}$.

\item If $a\leq0\leq b\leq a+3$, then  the functions $H_{a,b}$ are  convex in one direction and  univalent in $\mathbb{D}$.
Furthermore, if $a\leq0\leq b\leq a+2$, then  the functions $H_{a,b}$ are convex in $\mathbb{D}$.

\item The function $H_{a}$ is univalent in $\mathbb{D}$ if and only if $|a|\leq1$.
\end{enumerate}
\end{prop}
\bpf
(1) By computation, for all $z\in\mathbb{D}$, we have that
$$(1-|z|^2)|P_{H_{a,b}}(z)|=(1-|z|^2)\left |\frac{a+b+(b-a)z}{1-z^2}\right |\leq |a+b|+|b-a|=2\max\{|a|,|b|\}.
$$
Note that
$\lim_{r\rightarrow 1^-}(1-r^2)|P_{H_{a,b}}(r)|=2|b|$
and $\lim_{r\rightarrow -1^+}(1-r^2)|P_{H_{a,b}}(r)|=2|a|$.
Therefore, we get that
$\|P_{H_{a,b}}\|=2\max\{|a|,|b|\}$
and the result follows by Becker's univalence criterion.

Note that if $\max\{|a|,|b|\}>3$, then
$\|P_{H_{a,b}}\| >6$
and thus, the functions $H_{a,b}$ can not be univalent in $\mathbb{D}$.

(2) We observe that
$$H'_{a,b}(z)=\left(\frac{1+z}{1-z}\right)^a(1-z)^{a-b}
=\left(\frac{1+z}{1-z}\right)^b(1+z)^{a-b},
$$
and thus,
$$|\arg(H_{a,b}'(z))|<\frac{\pi}{2}\min\left\{(|a|+|a-b|),~ (|b|+|a-b|)\right\},\quad z\in\mathbb{D}.
$$
If $\min\{|a|,|b|\}+|a-b|\leq1$, then we have $|\arg(H_{a,b}'(z))|<\frac{\pi}{2}$ in $\mathbb{D}$ and thus,
by  Noshiro-Warschawski's theorem (see \cite{dur1983}), the functions
%which implies that $\RE H'_{a,b}(z)>0$ in $\mathbb{D}$. Clearly,
$H_{a,b}$ are close-to-convex and univalent in $\mathbb{D}$. %by Noshiro-Warschawski Theorem (see \cite{dur1983}).

(3) For $a\leq0\leq b\leq a+3$, %we will prove that $\RE (1+zH''_{a,b}(z)/H'_{a,b}(z))>-1/2$ in $\mathbb{D}$.
%(see \cite{oza} and \cite{rob}).
%called convex functions of order $-1/2$ and are univalent close-to-convex (see \cite{dur1} and \cite{goo}).
we see that
$$\RE \left(1+\frac{zH''_{a,b}(z)}{H'_{a,b}(z)}\right)
=1+\RE\frac{az}{1+z}+\RE\frac{bz}{1-z}>1+\frac{a}{2}-\frac{b}{2}=\frac{2+a-b}{2}\geq-\frac{1}{2}
$$
for all $z\in\mathbb{D}$ and thus, the functions $H_{a,b}$ are  convex in one direction and univalent in $\mathbb{D}$.
 Also, it is clearly that if $a\leq0\leq b\leq a+2$, then
$$\RE \left(1+\frac{zH''_{a,b}(z)}{H'_{a,b}(z)}\right)>0 ,\quad z\in\mathbb{D},
$$
and thus, the functions $H_{a,b}$ are convex in $\mathbb{D}$.

(4) It is a direct consequence of  \cite[Lemma~2.1]{KS2002} because of $H_a(z)=-H_{-a}(-z)$.
\epf

%  and $G_{a,b}(z)=H_{a+1,b}(z)-H_{a,b}(z)$.
%Obviously, the dilatation of $F_{a,b,\theta}$ is $e^{i\theta}z$.
%Thus, $F_{a,b,\theta}$ is sense-preserving harmonic mapping in $\mathbb{D}$ for all $a,~b,~\theta\in\mathbb{R}$.
%More precisely, $F_{a,b,\theta}\in \mathcal{A}_0$.
%For simplicity, if $a=b$, we denote $F_{a,a}$ and $G_{a,a}$ by $F_a$ and $G_a$ respectively. Thus, $F_a=H_a+\overline{G_a}$.

%\begin{lem}\label{LS3-lem2} Let $f = h + g$ be locally univalent in $\mathbb{D} $ and suppose that for
%some $\varepsilon$ ($|\varepsilon|\leq1$), $h + \varepsilon g$ is convex. Then $f$ is univalent  close-to-convex harmonic mapping in $D$. \end{lem}

%\begin{lem}\label{LS3-lem2}
%If $h$ and $g$ are analytic function defined on $\mathbb{D}$ with $h'(0)\neq0$,
% which satisfy $g'(z)=e^{i\theta}zh'(z)$ and $\RE \big(1+z\frac{h''}{h'}\big)>-\frac{1}{2}$ for all $z\in\mathbb{D}$,
% the harmonic function $f=h+\overline{g}$ is univalent and close-to-convex in $\mathbb{D}$.\end{lem}

\begin{prop}\label{LS3-prop5.2}
For all $\theta\in\mathbb{R}$, the family of harmonic mappings $F_{a,b,\theta}$ defined by \eqref{LS3-equ5.1}
has the following properties:
\begin{enumerate}
\item $||P_{F_{a,a,\theta}}||=2|a|+1=||P_{H_{a}}||+1$ for all $a\in\mathbb{R}$.

\item If $a\in(-\infty,-1]\cup[0,\infty)$, then $||P_{F_{a,a+1,\theta}}||=|2a+1|=||P_{H_{a,a+1}}||-1$.
         However, $||P_{F_{a,a+1,\theta}}||=1$ for each $a\in(-1,0)$.

\item If   $a\leq0\leq b\leq a+3$, then the functions $F_{a,b,\theta}$ are close-to-convex and univalent
           in $\mathbb{D}$.

\item   The functions $F_{a,a+1,\theta}$ are univalent in $\mathbb{D}$ if and only if   $a\in[-1,0]$.

%\item If   $a\in[-2,-1]$,  $F_{a,a+1,\frac{\pi}{2}}$ is univalent in $\mathbb{D}$.
\end{enumerate}
\end{prop}

\bpf By a straightforward computation, we have  that
$$P_{F_{a,b,\theta}}(z)
=\frac{a+b+(b-a)z}{1-z^2}-\frac{\overline{z}}{1-|z|^2},\quad z\in\mathbb{D}.
$$

(1) It follows from  \eqref{LS3-equ1.2} and \eqref{LS3-equ5.3} that $||P_{F_{a,a,\theta}}||=||P_{F_{-a,-a,\theta+\pi}}||$.
So we only need to consider the case  $a\geq0$.
The conclusion can be easily got by \eqref{LS3-equ2.1},
Proposition \ref{LS3-prop5.1} and the fact that
$$||P_{F_{a,a,\theta}}||\geq\lim_{r\rightarrow-1^+}(1-r^2)|P_{F_{a,a,\theta}}(r)|=2a+1=||P_{H_{a}}||+1,\quad a\geq0.
$$
%for $a\geq0$.

%$$
%\lim_{r\rightarrow1^-}|P_{F_{a,a,\theta}}(r)|(1-r^2)=|2a-1|=1-2a=1+||P_{H_{a,a}}||
%$$
%for $a<0$.
%Thus, we have $||P_{F_{a,a,\theta}}||\geq||P_{H_{a}}||+1$.
%%By Proposition \ref{LS3-prop5.1}, we know that $||P_{H_{a,a}}||=2|a|$.
%On the other hand, it follows from Theorem \ref{LS3-thm2.1} that %and Proposition \ref{LS3-prop5.1} that
%$||P_{F_{a,\theta}}||\leq||P_{H_{a}}||+1$.
%Thus, we get that $||P_{F_{a,\theta}}||=2a+1=||P_{H_{a}}||+1$.

(2) We first consider the case $a\geq0$.
Note that $$h_a(z):= H_{a,a+1}(z)+e^{i(\pi-\theta)}G_{a,a+1,\theta}(z)=H_a(z).
$$
%If $a\geq0$, then  $||P_{H_{a,a+1}}||=2(a+1)$ by Proposition \ref{LS3-prop5.1}.
It follows from \eqref{LS3-equ2.1} and Proposition \ref{LS3-prop5.1} that
$$2a+1=||P_{H_{a,a+1}}||-1\leq||P_{F_{a,a+1,\theta}}||\leq||P_{h_a}||+1=2a+1,\quad a\geq0.
$$
Obviously, $||P_{F_{a,a+1,\theta}}||=2a+1=||P_{H_{a,a+1}}||-1$ for each $a\geq0$ and all $\theta\in\mathbb{R}$.
For the case $a\leq-1$, the conclusion follows, since $||P_{F_{a,a+1,\theta}}||=||P_{F_{-(a+1),-(a+1)+1,\theta+\pi}}||$.
%$||P_{a,a+1,\theta}||=||P_{-(a+1),-(a+1)+1,\theta+\pi}||$.
%by \eqref{LS3-equ1.2} and \eqref{LS3-equ5.3}.
%To acquire the equality $||P_{F_{a,a+1,\theta}}||=2a+1$,
%it suffices to prove %$||P_{F_{a,a+1,\theta}}||\leq2a+1$
% $|P_{F_{a,a+1,\theta}}|(1-|z|^2)\leq2a+1$,
%i.e.,
%$$
%\left|\frac{(2a+1)(1-|z|^2)+z-\overline{z}}{1-z^2}\right|\leq2a+1,
%$$
%which can be verified from the following inequality
%\begin{align*}
%& &(2a+1)^2|1-z^2|^2-((2a+1)^2(1-|z|^2)^2+(z-\overline{z})^2)\\
%&=&2((2a+1)^2(|z|^2- \RE z^2))-2\IM^2z)\\
%&=& 4((2a+1)^2-1)\IM^2z\geq0.
%\end{align*}
%By the similar analysis as above, we can get the corresponding conclusion for $a\leq-1$.

Next we will certify that $||P_{F_{a,a+1,\theta}}|| =1$ for each $a\in(-1,0)$.
A basic computation states that
$$|1-z^2|^2-((2a+1)^2(1-|z|^2)^2-(z-\overline{z})^2)=-4a(a+1)(1-|z|^2)^2>0,\quad z\in\mathbb{D},
$$
%&=&-4a(a+1)(1+|z|^4)+2((2a+1)^2|z|^2-\RE z^2)\\
%&=&-4a(a+1)(1-x^2)^2+((2a+1)^2+1-4a(a+1)(2x^2+y^2))y^2\geq0,
which means that
$$(1-|z|^2)|P_{F_{a,a+1,\theta}}(z)|=\Big|\frac{(2a+1)(1-|z|^2)+z-\overline{z}}{1-z^2}\Big|\leq1,\quad z\in\mathbb{D}.
$$
It yields that $||P_{F_{a,a+1,\theta}}||\leq1$.
Note that $\lim_{r\rightarrow1^-}(1-r^2)|P_{F_{a,a+1,\theta}}(ir)|=1$ and thus,
we obtain $||P_{F_{a,a+1,\theta}}||=1$.

(3)  If   $a\leq0\leq b\leq a+3$, then from the proof of Proposition \ref{LS3-prop5.1} we find that
$$\RE \left(1+\frac{zH''_{a,b}(z)}{H'_{a,b}(z)}\right)>-\frac{1}{2}, \quad z\in\mathbb{D}.
$$
%$\RE (1+zH''_{a,b}(z)/H'_{a,b}(z))>-1/2$ in $\mathbb{D}$
Note that the dilatation of $F_{a,b,\theta}$ is $e^{i\theta}z$ for all $a,b\in\mathbb{R}$ and each $\theta\in\mathbb{R}$.
As a consequence, it follows from \cite[Theorem~1]{BP2014} that the functions $F_{a,b,\theta}$
are close-to-convex and univalent for all $\theta\in\mathbb{R}$ if $a\leq0\leq b\leq a+3$.

(4)  Obviously, it follows from (3) that the functions $F_{a,a+1,\theta}$ are univalent in $\mathbb{D}$
for all $\theta\in\mathbb{R}$  if $a\in[-1,0]$. For the remaining part, we use the method of contradiction.
Assume that there exists some  $a>0$ such that the functions $F_{a,a+1,\theta}$ are univalent in $\mathbb{D}$
for all $\theta\in\mathbb{R}$.
Therefore, the function $F_{a,a+1,\theta}$ is stable univalent for each $\theta\in\mathbb{R}$
and thus $H_{a,a+1}+\lambda G_{a,a+1,\theta}$ is univalent in $\mathbb{D}$ for each $\lambda\in\overline{\mathbb{D}}$ (see \cite{HM2013}),
especially for $\lambda=e^{-i\theta}$.
However, $H_{a,a+1}+e^{-i\theta}G_{a,a+1,\theta}=H_{a+1,a+1}$ is not univalent in $\mathbb{D}$ by Proposition \ref{LS3-prop5.1} when $a>0$.
This is a contradiction.
Using \eqref{LS3-equ5.3}, the similar contradiction can be obtained  for the case $a<-1$.
This completes the proof.
%However, for $a>0$, we have that
%\begin{align*}
%||S_{F_{a,a+1,0}}||=&\sup_{z\in\mathbb{D}}(1-|z|^2)^2|S_{F_{a,a+1,0}}(z)|\\
%=&\sup_{z\in\mathbb{D}}(1-|z|^2)^2\frac{2(a+1)|2z-a-1|}{|1-z^2|^2}=2(a+1)(a+3)>6,
%\end{align*}
%which implies that $F_{a,a+1,0}$ is not univalent in $\mathbb{D}$. % by Nehari-Kraus's criterion.
%Again, for $a<-1$,  it follows from \eqref{LS3-equ5.3} that
%$$||S_{F_{a,a+1,\pi}}||=||S_{F_{-(a+1),-(a+1)+1,2\pi}}||=2|a|(|a|+2)>6,
%$$
%which gives that $F_{a,a+1,\pi}$ is not univalent in $\mathbb{D}$. % because of  $||S_{F_{a,a+1,\pi}}||=2|a|(|a|+2)>6$.
 \epf

From the proof of Proposition \ref{LS3-prop5.2},
the two families of harmonic mappings $F_{a,a,\theta}$ and $F_{a,a+1,\theta}$
provide sharp results for several of the inequalities in Section \ref{LS3-sec2}.
For simplicity, let
$$h_{a,b,\theta,\varphi}(z)=H_{a,b}(z)+e^{i\varphi}G_{a,b,\theta}(z)
=H_{a,b}(z)+e^{i(\theta+\varphi)}(H_{a+1,b}(z)-H_{a,b}(z)).%\quad z\in\mathbb{D}.
$$
Clearly, $h_{a,b,\theta,-\theta}=H_{a+1,b}$
and $h_{a,b,\theta,\pi-\theta}=2H_{a,b}-H_{a+1,b}=H_{a,b-1}$.
%Based on the conclusions in Proposition \ref{LS3-prop5.1} and \ref{LS3-prop5.2}
%about pre-Schwarzian norms of $H_{a,b}$ and $F_{a,b,\theta}$,
We have the following results from Propositions \ref{LS3-prop5.1} and \ref{LS3-prop5.2}.
\bei
\item $||P_{h_{a,a,\theta,\pi-\theta}}||+1=||P_{H_a}||+1=||P_{F_{a,a,\theta}}||
=2a+1=||P_{h_{a,a,\theta,-\theta}}||-1$,~~$a\geq1/2$.\\
\item $||P_{H_a}||+1=||P_{F_{a,a,\theta}}||
=2a+1=||P_{h_{a,a,\theta,-\theta}}||-1$,~~$0\leq a\leq1/2$.\\
\item $||P_{F_{a,a,\theta}}||=2a+1=||P_{h_{a,a,\theta,\pi-\theta}}||+\varepsilon$,
~~$a=(1+\varepsilon)/4\in[1/4,1/2]$.\\
\item $||P_{F_{a,a,\theta}}||=2a+1=||P_{h_{a,a,\theta,\pi-\theta}}||-\varepsilon$,
~~$a=(1-\varepsilon)/4\in[0,1/4]$.\\
\item $||P_{h_{a,a+1,\theta,\pi-\theta}}||+1=||P_{F_{a,a+1,\theta}}||
=2a+1=||P_{H_{a,a+1}}||-1=||P_{h_{a,a+1,\theta,-\theta}}||-1$,~~$a\geq0$.\\
\eei
Similar results may be stated for $F_{a,a,\theta}$  and
$F_{a,a+1,\theta}$ when $a<0$. %by \ref{LS3-equ5.3}.
For these functions, we know that $||P_{F_{a,b,\theta}}||\geq1$
and $\omega_{F_{a,b,\theta}}=e^{i\theta}z$ for all $a$, $b=a~\text{or}~a+1$, $\theta\in\mathbb{R}$.
These things do not happen accidentally.
Our next result, which is a parallel result to \cite[Theorem~3]{CHM}, demonstrates the reason behind these.
%\begin{cor} Furthermore,  for any given $\lambda\geq1$, there exist some $f_\lambda=h_\lambda+\overline{g_\lambda}\in\mathcal{H}^0_\lambda$ and some $\varepsilon\in\overline{\mathbb{D}}$ such that $||P_{h_\lambda+\varepsilon g_\lambda}||=||P_{f_\lambda}||+1$ or $||P_{f_\lambda}||=||P_{h_\lambda+\varepsilon g_\lambda}||+1$.\end{cor}

%If $f\in\mathcal{H}^0_\lambda$, the following proposition will prove that the constant 1 in \eqref{LS3-eq1} can not be got.

%Set $\mathcal{B}_H^*(\lambda)=
%\{f\in\mathcal{B}_H(\lambda)~~\text{with}~~||P_f||=\lambda\}$
%and $\mathcal{B}_{H}^{0*}(\lambda)=\mathcal{B}_H^*(\lambda)\cap\mathcal{H}_0$.
We denote by $\mathcal{A}(\lambda)$ (resp. $\mathcal{A}_0(\lambda)$) the set of all {\it admissible dilatations}
of $f\in\mathcal{B}_H(\lambda)$ (resp. $\mathcal{B}_{H}^0(\lambda))$;
i.e., $\omega\in\mathcal{A}(\lambda)$ (or $\mathcal{A}_0(\lambda)$) if there exists a harmonic mapping
$f=h+\overline{g}\in\mathcal{B}_H(\lambda)$
$(\mathcal{B}_{H}^0(\lambda))$ with dilatation $\omega$.
%The following proposition characterizes the values of $\lambda$ for which a dilatation in $\mathcal{A}^*_\lambda$ can have hyperbolic norm 1.

\bthm\label{LS3-thm5.1}
The following conditions are equivalent.
\begin{enumerate}
\item $\lambda\geq1$;

\item There exists a $\omega\in\mathcal{A}_0(\lambda)$ with $|\omega'(0)|=1$;

\item The set $\{\mu\cdot I:\,|\mu|=1\}$ is contained in $\mathcal{A}_0(\lambda)$;
%In particular, the identity function $I$ is an admissible dilatation in $\mathcal{B}^0_{\lambda}$;

\item Every automorphism $\sigma$ of the unit disk is an admissible dilatation in
$\mathcal{B}_H(\lambda)$.
\end{enumerate}
\ethm
\bpf
The scheme of the proof is to show that $(1)\Rightarrow(2)\Rightarrow(3)\Rightarrow(1)$
and $(3)\Leftrightarrow(4)$.
We only show $(1)\Rightarrow(2)$ and $(3)\Rightarrow(1)$.
The remaining implications are similar to corresponding proofs of \cite[Theorem~3]{CHM}
since the PSD preserves affine invariance.

We now show that $(1)\Rightarrow(2)$: For any given $\lambda\geq1$, %for any given $\lambda\geq1$,
we  choose $|a|=\frac{\lambda-1}{2}$ so that $P_{F_{a,a,\theta}}\in\mathcal{B}_{H}^0(\lambda)$
with the dilatation $e^{i\theta}z$ by Proposition \ref{LS3-prop5.2}.  Then (2) follows.
%Obviously, $|(e^{i\theta}z)'(0)|\equiv1$.

Next, we prove that $(3)\Rightarrow(1)$:
If $(3)$ is satisfied, then for any given $\mu$ with $|\mu|=1$,
there is a harmonic function $f_{\mu}=h_{\mu}+\overline{g_{\mu}}
\in\mathcal{B}_{H}^0(\lambda)$ with dilatation $\mu z$.
Since $g'_\mu(0)=0$, by \eqref{LS3-equ1.3},
  %$||P_{h_{\mu}+\varepsilon g_{\mu}}||\leq||P_{f_{\mu}}||+1=\lambda+1$
% for any $\varepsilon\in\overline{\mathbb{D}}$.
%Note that $(h_{\mu}+\varepsilon g_{\mu})(0)=(h_{\mu}+\varepsilon g_{\mu})'(0)-1=0$.
 $h_{\mu}+\varepsilon g_{\mu}\in\mathcal{B}_{A}((\lambda+1)/2)$
for each $\varepsilon\in\overline{\mathbb{D}}$.
It follows from  \cite[Theorem~2.3]{KS2002} that
for any $\varepsilon\in\overline{\mathbb{D}}$,
\be \label{LS3-equ5.4}
|(h_{\mu}'+\varepsilon g_{\mu}')(z)|=|h_{\mu}'(z)|\cdot|1+\varepsilon\mu z|
\geq\left(\frac{1-|z|}{1+|z|}\right)^{\frac{\lambda+1}{2}},\quad z\in\mathbb{D}.
\ee
Since $\varepsilon\in\overline{\mathbb{D}}$ and $|\mu|=1$, for any given $z\neq0$ in the unit disk, we can get
$$|h_{\mu}'(z)|\geq\frac{(1-|z|)^{\frac{\lambda-1}{2}}}{(1+|z|)^{\frac{\lambda+1}{2}}},
$$
by choosing $\varepsilon\cdot\mu=-\overline{z}/|z|$ in \eqref{LS3-equ5.4}.
Clearly, the above inequality holds for $z=0$.
Note that $h_{\mu}$ is locally univalent in $\mathbb{D}$
and $h'_{\mu}(0)=1$ for each $\mu\in\partial\mathbb{D}$.
We obtain that the analytic function $1/h'_{\mu}$ satisfies
$$\frac{1}{|h_{\mu}'(z)|}\leq\frac{(1+|z|)^{\frac{\lambda+1}{2}}}{(1-|z|)^{\frac{\lambda-1}{2}}},\quad z\in\mathbb{D},
$$
which implies that $\lambda\geq1$ by the maximum modulus principle.
Otherwise we would get $1/|h_{\mu}'(0)|<1$, which contradicts  $h'_{\mu}(0)=1$.
This completes the proof.
\epf

Compared to \cite[Theorem~3]{CHM},
since the PSD is in general not linear invariant,
we are not sure whether the conditions in Theorem \ref{LS3-thm5.1}
are equivalent to that there exists $\omega\in\mathcal{A}(\lambda)$ (or $\mathcal{A}_0(\lambda)$) with $||\omega^*||=1$.
 Here $||\omega^*||$ is the hyperbolic norm of the dilatation $\omega$ of a sense-preserving harmonic mapping in $\mathbb{D}$, i.e.,
$$||\omega^*||=\sup_{z\in\mathbb{D}}\frac{|\omega'(z)|(1-|z|^2)}{1-|\omega(z)|^2}.
$$
The hyperbolic norm plays a distinguished role in the analysis of the order of affine and linear invariant families
of harmonic mappings with bounded SD (see \cite{CHM}).

\section{Growth estimate for the class $\mathcal{B}_{H}(\lambda)$} \label{LS3-sec6}
To study the growth estimate for the class $\mathcal{B}_{H}(\lambda)$,
we need the following result which characterizes harmonic mappings in $\mathcal{B}_{H}(\lambda)$.

\begin{prop}\label{LS3-prop6.1}
A  harmonic mapping $f \in\mathcal{H}$ belongs to $\mathcal{B}_{H}(\lambda)$
if and only if for each pair of points $z$, $z_0$ in $\mathbb{D}$, the inequality
$$|A(z)-A(z_0)|\leq\lambda d_h(z,z_0)
$$
holds, where $A(z)=\log J_f(z)$.
\end{prop}
\bpf %Let $f=h+\overline{g}\in\mathcal{H}$ be a locally univalent harmonic mapping  and
Assume that  $f=h+\overline{g}\in\mathcal{B}_{H}(\lambda)$.
Then $|P_f(z)|\leq \lambda /(1-|z|^2)$ holds in $\mathbb{D}$.
We observe that
$$A_{\overline{z}}=(\log J_f)_{\overline{z}} =\overline{(\log J_f)_z}=\overline{P_f}=\overline{A_z}.
$$
Therefore, for two points $z$, $z_0$ in $\mathbb{D}$, we have that
\begin{align*}
 |A(z)-A(z_0)|
\leq&\left|\int_\Gamma A_\zeta(\zeta)d\zeta
+A_{\overline{\zeta}}(\zeta)d\overline{\zeta}\right|\\
\leq&\int_\Gamma(|A_\zeta(\zeta)|
+|A_{\overline{\zeta}}(\zeta)|)|d\zeta|\\
\leq&\int_\Gamma\frac{2\lambda}{1-|\zeta|^2}|d\zeta|
=\lambda d_h(z,z_0),
\end{align*}
where $\Gamma$ is the hyperbolic geodesic joining $z$ and $z_0$.

Conversely, we assume that the inequality
$|A(z)-A(z_0)|\leq\lambda d_h(z_1,z_2)$
holds for each pair of points $z$, $z_0$ in $\mathbb{D}$.
%In order to get that $||P_f||\leq\lambda$,
It suffices to prove that
 $(1-|z_0|^2)|P_f(z_0)|\leq\lambda$ for each $z_0\in\mathbb{D}$.
Fix $z_0\in\mathbb{D}$. If  $A_z(z_0)=0$, then $P_f(z_0)=0$ and $(1-|z_0|^2)|P_f(z_0)|\leq\lambda$.
Otherwise, choose a curve
$$\gamma =\{z=z_0+re^{-i\theta}:\, r\in(0,1-|z_0|),~\theta=\arg(A_z(z_0))\}.
$$
Clearly, $A(z)$ is infinitely differentiable in $\mathbb{D}$ owing to $J_f(z)>0$. %and $h$, $g$ are analytic in $\mathbb{D}$.
Thus, we have the representation
$$A(z)-A(z_0)=A_z(z_0)(z-z_0)+A_{\overline{z}}(z_0)(\overline{z}-\overline{z_0})
+\sum_{i+j>1}C_{ij}(z-z_0)^i(\overline{z}-\overline{z_0})^j
$$
for some complex constants $C_{ij}$, which implies that
$$\lim_{\gamma\ni z\rightarrow z_0}\frac{A(z)-A(z_0)}{z-z_0}
%=A_z(z_0)+A_{\overline{z}}(z_0)e^{2i\theta}\\
=A_z(z_0)+\overline{A_z(z_0)}e^{2i\theta}=2A_z(z_0)=2P_f(z_0).
$$
The desired inequality $(1-|z_0|^2)|P_f(z_0)|\leq\lambda$ follows from the equality
$$\lim_{\gamma\ni z\rightarrow z_0}\frac{|A(z)-A(z_0)|}{d_h(z,z_0)}=(1-|z_0|^2)|P_f(z_0)|.
$$
The proof is complete.
\epf

\bthm\label{LS3-thm6.1} {\rm (Distortion theorem)} Let
$f=h+\overline{g}\in\mathcal{B}_{H}(\lambda)$ for some $\lambda\geq0$ with $b_1=g'(0)$, and let $H_{a,b}$ and $H_a$ be
defined by \eqref{LS3-equ5.2}. Then for each $z\in\mathbb{D}$, we have
\begin{enumerate}
\item $(1-|b_1|^2)H'_\lambda(-|z|)\leq J_f(z)\leq(1-|b_1|^2)H'_\lambda(|z|)$;

%$(1-|b_1|^2)H'_\lambda(-|z|)=\big(\frac{1-|z|}{1+|z|}\big)^\lambda\leq J_f(z)\leq\big(\frac{1+|z|}{1-|z|}\big)^\lambda=H'_\lambda(|z|)$;

\item $\sqrt{1-|b_1|^2}H'_\frac{\lambda}{2}(-|z|)\leq |h'(z)|
        \leq (1+|b_1z|)H'_{\frac{\lambda-1}{2},\frac{\lambda+1}{2}}(|z|)$;

\item $|g'(z)|\leq (|z|+|b_1|)H'_{\frac{\lambda-1}{2},\frac{\lambda+1}{2}}(|z|)$;

\item $-\sqrt{1-|b_1|^2}H_\frac{\lambda}{2}(-|z|)
        \leq |h(z)|\leq (1-|b_1|)H_{\frac{\lambda-1}{2},\frac{\lambda+1}{2}}(|z|)+|b_1|H_{\frac{\lambda+1}{2}}(|z|)$;

\item $|g(z)|\leq H_\frac{\lambda+1}{2}(|z|)-(1-|b_1|)H_{\frac{\lambda-1}{2},\frac{\lambda+1}{2}}(|z|)$;

\item $|f(z)|\leq (1+|b_1|)H_{\frac{\lambda+1}{2}}(|z|)$.
Furthermore, if $f$ is univalent in $\mathbb{D}$, then
$$-(1-|b_1|)H_{\frac{\lambda+1}{2}}(-|z|)
\leq|f(z)|\leq (1+|b_1|)H_{\frac{\lambda+1}{2}}(|z|).
$$
\end{enumerate}

The estimates in {\rm (1)} are sharp for all $\lambda\geq0$. The right sides of {\rm (2)-(6)} are sharp for
all $\lambda\geq1$ and the left side of {\rm (6)} is sharp for $\lambda=1$.
Moreover, if $f\in\mathcal{B}_H^0(\lambda)$,
then the left sides of {\rm (2)} and {\rm (4)} are sharp for all  $\lambda\geq0$.
\ethm

\bpf
(1) The conclusion can be easily obtained by choosing $z_0=0$ in Proposition \ref{LS3-prop6.1}.

(2) Since $f\in\mathcal{B}_H(\lambda)$, %then $|\omega_f|<1$ and $\omega(0)=0$.
by Lindel\"{o}f's inequality, we get that
$$|\omega_f(z)|\leq\frac{|z|+|b_1|}{1+|b_1z|}
$$
and thus,
\begin{align*} |h'(z)|&=\left(\frac{J_f(z)}{1-|w_f(z)|^2}\right )^{\frac{1}{2}}\\
&\leq (1-|b_1|^2)^\frac{1}{2}\left (1-\left (\frac{|z|+|b_1|}{1+|b_1z|}\right )^2\right )^{-\frac{1}{2}}
\left (\frac{1+|z|}{1-|z|}\right )^{\frac{\lambda}{2}}\\
&=(1+|b_1z|)H'_{\frac{\lambda-1}{2},\frac{\lambda+1}{2}}(|z|)
%&=(1+|b_1z|)\frac{(1+|z|)^{\frac{\lambda-1}{2}}}{(1-|z|)^{\frac{\lambda+1}{2}}}
%\leq(1+|b_1|)\frac{(1+|z|)^{\frac{\lambda-1}{2}}}{(1-|z|)^{\frac{\lambda+1}{2}}}
%=(1+|b_1|)H'_{\frac{\lambda-1}{2},\frac{\lambda+1}{2}}(|z|)
\end{align*}
and
$$|h'(z)|\geq(J_f(z))^\frac{1}{2}\geq\sqrt{1-|b_1|^2}\left (\frac{1-|z|}{1+|z|}\right )^{\frac{\lambda}{2}}
=\sqrt{1-|b_1|^2}H'_{\frac{\lambda}{2}}(-|z|).
$$
%If $f\in\mathcal{B}_H^0(\lambda)$ ($\lambda\geq1)$, the sharp example can be seen from
%$F_{\frac{\lambda-1}{2},\frac{\lambda+1}{2},\theta}$ for all $\theta\in\mathbb{R}$
%and $H_{\frac{\lambda}{2}}$  in turn since $b_1=0$.

(3) It follows from Lindel\"{o}f's  inequality and the proof of (2) that
$$|g'(z)|=|w_f(z)h'(z)|%\leq(|z|+|b_1|)\frac{(1+|z|)^{\frac{\lambda-1}{2}}}{(1-|z|)^{\frac{\lambda+1}{2}}}
\leq(|z|+|b_1|)H'_{\frac{\lambda-1}{2},\frac{\lambda+1}{2}}(|z|).
$$
%If $f\in\mathcal{B}_H^0(\lambda)$ ($\lambda\geq1)$, the estimate is sharp by choosing
%$f=F_{\frac{\lambda-1}{2},\frac{\lambda+1}{2},\theta}$ for all $\theta\in\mathbb{R}$.

(4) Integrating inequalities in (2) yields (4).

(5) Integrating inequality in (3) yields (5).

(6) Applying the triangle inequality and the results in (4) and (5), we obtain
$$|f(z)|\leq|h(z)|+|g(z)|\leq (1+|b_1|)H_{\frac{\lambda+1}{2}}(z).
$$

Let $f_\varepsilon=\frac{h+\varepsilon g}{1+\varepsilon b_1}~(\varepsilon\in\overline{\mathbb{D}})$. Then $f_\varepsilon$
belongs to $\mathcal{B}_{A}(\frac{\lambda+1}{2})$.
By \cite[ Theorem~2.3]{KS2002}, we have % that for any $\varepsilon\in\overline{\mathbb{D}}$ and $z\in\mathbb{D}$,
$$|(h'+\varepsilon g')(z)|\geq|1+\varepsilon b_1|\left (\frac{1-|z|}{1+|z|}\right )^{\frac{\lambda+1}{2}}
\geq%(1-|b_1|)\Big(\frac{1-|z|}{1+|z|}\Big)^{\frac{\lambda+1}{2}}=
(1-|b_1|)H'_{\frac{\lambda+1}{2}}(-|z|).
$$
Especially, since $\varepsilon$ is arbitrary, we get that
$$|h'(z)|-|g'(z)|\geq(1-|b_1|)H'_{\frac{\lambda+1}{2}}(-|z|).%,\quad z\in\mathbb{D}.
$$
For $0<r<1$ we choose $z_0$ such that $|f(z_0)|$ is the minimum of $|f(z)|$ on $|z|=r$.
If $f$ is univalent in $\mathbb{D}$ and  $\gamma$ is the preimage of the segment $[0,f(z_0)]$, then for $|z|=r$, we have that
$$|f(z)|\geq|f(z_0)|=\int_{\gamma}|df(z)| \geq\int_0^{|z|}(|f_z(z)|
-|f_{\overline{z}}(z)|)|dz|\geq-(1-|b_1|)H_{\frac{\lambda+1}{2}}(-|z|).
$$
%\begin{align*} |f(z)|&\geq|f(z_0)|\geq\int_{f(\Gamma))}|df(z)|\geq\int_0^{|z|}(|f_z(z)|-|f_{\overline{z}}(z)|)|dz|\geq-(1-|b_1|)H_{\frac{\lambda+1}{2}}(-|z|)\\
%&\geq\int_0^{|z|}\Big(\frac{1-t}{1+t}\Big)^{\frac{\lambda+1}{2}}dt=-(1-|b_1|)H_{\frac{\lambda+1}{2}}(-|z|).
%\end{align*}

Next we consider the sharpness part.
The equality occurs in (1) if we take
$$f(z)=H_{\frac{\lambda}{2}}(z)+\overline{b_1H_{\frac{\lambda}{2}}(z)} ~\mbox{ for each $\lambda\geq0$}.
$$
For each $\lambda\geq1$, the equalities in the right sides of (2)-(6) are attained for
$$f(z)=f_\lambda (z)=F_{\frac{\lambda-1}{2},\frac{\lambda+1}{2},0}(z)+|b_1|\overline{F_{\frac{\lambda-1}{2},\frac{\lambda+1}{2},0}(z)}
$$
at $z=r\in[0,1)$, where $F_{a,b,0}$ is defined by \eqref{LS3-equ5.1}.
Note that $f_\lambda\in\mathcal{B}_H(\lambda)$ by \eqref{LS3-equ1.2} and Proposition \ref{LS3-prop5.2}.
Similarly, for each $\lambda\geq0$, the function $H_{\frac{\lambda}{2}}$ provides the sharpness
for the left sides of (2) and (4) at $z=-r\in(-1,0]$ when $f\in\mathcal{B}_H^0(\lambda)$.
It follows from Proposition \ref{LS3-prop5.2} that $F_{\frac{\lambda-1}{2},\frac{\lambda+1}{2},0}$
is univalent in $\mathbb{D}$ for $\lambda=1$.
The equality in the left side of (6) occurs for
%$f=F_{\frac{\lambda-1}{2},\frac{\lambda+1}{2},0}\in\mathcal{B}_H^0(\lambda)$
$f=F_{0,1,0}-|b_1|\overline{F_{0,1,0}}\in\mathcal{B}_H(1)$ and  $z=-r\in(-1,0]$.
We complete the proof.
\epf

The following result can be directly deduced from Theorem \ref{LS3-thm6.1}.

\bcor \label{LS3-cor6.1} {\rm (Growth and covering theorem)}
Let $f=h+\overline{g}\in\mathcal{B}_{H}(\lambda)$ with $b_1=g'(0)$, and let
$H_{a,b}$ and $H_a$ be defined by \eqref{LS3-equ5.2}.
If $\lambda>1$, then $f$, $h$ and $g$ satisfy the same growth condition
$$f(z)~(h(z),~g(z))=O(1-|z|)^{\frac{1-\lambda}{2}} ~\mbox{ as $|z|\rightarrow 1$.}
$$
If $\lambda<1$, then $f$ (resp. $h$, $g$) is bounded by
$$(1+|b_1|)H_{\frac{\lambda+1}{2}}(1) ~(resp. ~(1-|b_1|)H_{\frac{\lambda-1}{2},\frac{\lambda+1}{2}}(1)+|b_1|H_{\frac{\lambda+1}{2}}(1),
~H_{\frac{\lambda+1}{2}}(1)-(1-|b_1|)H_{\frac{\lambda-1}{2},\frac{\lambda+1}{2}}(1)).
$$

For all $\lambda>0$, the image  $h(\mathbb{D})$ contains the disk $\{|z|<-\sqrt{1-|b_1|^2}H_{\frac{\lambda}{2}}(-1)\}$.
If $f\in\mathcal{B}_{H}(\lambda)\cap\mathcal{S}_{H}$,
then the image $f(\mathbb{D})$  contains the disk $\{|z|<-(1-|b_1|)H_{\frac{\lambda+1}{2}}(-1)\}$.
\ecor

%We also note that, for $0\leq \lambda\leq 1$, we have
%The former part and the conclusion of the image $h(\mathbb{D})$ can be directly deduced from Theorem \ref{LS3-thm6.1}.
%For the image $f(\mathbb{D})$, the proof  is similar to (6) in Theorem \ref{LS3-thm6.1}.

If $f\in\mathcal{B}_{H}(\lambda)\cap\mathcal{S}^0_{H}$ for some $\lambda\in[0,1]$,
then
$$-H_{\frac{\lambda+1}{2}}(-1)\geq-H_1(-1)=2\log 2-1=0.38629\cdots.
$$
%and $-H_{\frac{\lambda}{2}}(-1)\geq-H_{\frac{1}{2}}(-1)=\frac{\pi}{2}-1=0.57082\cdots$.
This result is an improvement over the non-sharp known result that $f(\mathbb{D})\supseteq\{w:\,|w|<1/16\}$ if $f\in\mathcal{S}_H^0$.

%of the covering theorem for harmonic mappings in $\mathcal{S}_H^0$
%in the sense that if $f\in\mathcal{S}_H^0$, then the non-sharp known result is that $f(\mathbb{D})\supseteq\{w:\,|w|<1/16\}$.

In Corollary \ref{LS3-cor6.1}, the case $\lambda=1$ is critical.
By Theorem \ref{LS3-thm6.1}, we have that, for $f\in \mathcal{B}_H(1)$,
$$|f(z)|\leq (1+|b_1|)H_1(|z|)=(1+|b_1|)\big(-2\log(1-|z|)-|z|\big),\quad z\in\mathbb{D},
$$
which shows that functions in $\mathcal{B}_H(1)$ need not be bounded.
The following result gives a sufficient condition for the boundedness  of mappings in $\mathcal{B}_H(1)$.

\begin{prop} \label{LS3-prop6.2}
Let $f=h+\overline{g}$ be a sense-preserving harmonic mapping in $\mathbb{D}$.
If $f$  satisfies the condition
$$\beta(f):=\overline{\lim_{|z|\rightarrow1^{-}}}((1-|z|^2)|P_f(z)|-1)\log\frac{1}{1-|z|^2}<-2,
$$
then $f$, $h$ and $g$ are bounded in $\mathbb{D}$.
%Here the constant -2 on the right side is sharp.
\end{prop}
\bpf
Without loss of generality, we can assume that $g(0)=0$.
It follows from \eqref{LS3-equ1.2} that $\beta(A\circ f)=\beta(f)$ for any affine harmonic mapping
$A$ defined in \eqref{LS3-equ1.2}. Let $A\circ f=H+\overline{G}$.
It is easy to check that both $h$ and $g$ are bounded in $\mathbb{D}$ if and only if
both $H$ and $G$ are bounded in $\mathbb{D}$.
Note that $A\circ f$ is also sense-preserving in $\mathbb{D}$.
%Based on those with the analysis in section \ref{LS3-sec1} and hypothesis,
Thus, it is enough to consider the case $f=h+\overline{g}\in\mathcal{H}_0$ %\mathcal{H}^0$
and prove that both $h$ and $g$ are bounded in $\mathbb{D}$.

By assumption, there exist $\beta<-2$ and $r_0\in(1-1/(2e),1)$ such that
\be \label{LS3-equ6.1}
|P_f(z)|\leq\frac{1}{1-|z|^2}+\frac{\beta}{(1-|z|^2)\log(1/1-|z|^2)}
\ee
for $z\in D_{r_0}=\{z:\,r_0<|z|<1\}$.
Fix $z\in D_{r_0}$ and let $\Gamma$ be a line segment from $z$ to $z_0:=r_0e^{i\arg z}$ in the proof of
Proposition \ref{LS3-prop6.1}. Then we have
\be \label{LS3-equ6.2}
|\log J_f(z)|\leq2\int_{r_0}^{|z|}|P_f(\zeta)||d\zeta|+C_1,
\ee
where $C_1=\max_{\theta\in[0,2\pi]}|J_f(r_0e^{i\theta})|<\infty$.
By (\ref{LS3-equ6.1}) and (\ref{LS3-equ6.2}), we see that
\begin{align*}
 |\log J_f(z)|&\leq\log\frac{1+|z|}{1-|z|}+\int_{r_0}^{|z|}\frac{2\beta dt}{(1-t^2)\log(1/(1-t^2))}+C_1\\
&\leq\log\frac{1+|z|}{1-|z|}+\int_{r_0}^{|z|}\frac{\beta dt}{(1-t)\log(1/(2(1-t)))}+C_1\\
&=\log\frac{1+|z|}{1-|z|}+\beta\log\log\frac{1}{2(1-|z|)}+C_2,
\end{align*}
where $C_2=C_1-\beta\log\log\frac{1}{2(1-r_0)}$. % <\infty$.
Exponentiating the last inequality shows that
$$|J_f(z)|=|h'(z)|^2(1-|\omega_f(z)|^2)\leq e^{C_2}\frac{1+|z|}{1-|z|}\left(\log\frac{1}{2(1-|z|)}\right)^\beta.
$$
Using that $f\in\mathcal{B}_H^0$(1), we have $|\omega_f(z)|\leq|z|$ in $\mathbb{D}$ and thus, we find that
$$|g'(z)|< |h'(z)|=\left(\frac{J_f(z)}{1-|\omega_f(z)|^2}\right)^{1/2}\leq
 \frac{e^{C_2/2}}{1-|z|}\left(\log\frac{1}{2(1-|z|)}\right)^{\beta/2}.
$$
Since $\beta/2<-1$, the function $[\log(1/(2(1-t)))]^{\beta/2}/(1-t)$ is integrable on the interval $[r_0,1)$.
It follows that both $h$ and $g$ are bounded in $\mathbb{D}$ so that $f$ is also bounded in $\mathbb{D}$.
\epf

\br
Let $f=h+\overline{g}$ be a sense-preserving harmonic mapping in $\mathbb{D}$.
If $h$ and $g$ are unbounded in $\mathbb{D}$, then the boundedness of $f$ is uncertain.
For instance, let's recall the function $F_{0,1,\theta}$ defined by \eqref{LS3-equ5.1}:
$$F_{0,1,\theta}(z)=H_{0,1}(z)+\overline{G_{0,1,\theta}(z)}=-\log(1-z)+\overline{e^{i\theta}(-z-\log(1-z))},
$$
we see that $H_{0,1}$ and $G_{0,1,\theta}$ are unbounded in  $\mathbb{D}$.
However, $F_{0,1,0}(z)=-\overline{z}-2\log|1-z|$ and $F_{0,1,\pi}(z)=\overline{z}-2\arg(1-z)$
are unbounded and bounded in  $\mathbb{D}$, respectively.
Furthermore, it follows from Proposition \ref{LS3-prop5.2} that $||P_{F_{0,1,\theta}}||=1$ for any $\theta\in\mathbb{R}$.
\er

By Theorem \ref{LS3-thm6.1} and \cite[Theorem~5.1]{dur1970},  we conclude
 the H\"{o}lder continuity of mappings in $\mathcal{B}_H(\lambda)$.

\bthm\label{LS3-thm6.3}
Let $f=h+\overline{g}\in\mathcal{B}_H(\lambda)$ for some $\lambda\in[0,1)$.
Then $h+\varepsilon g$ is H\"{o}lder continuous of exponent $\frac{1-\lambda}{2}$ in $\mathbb{D}$ for each
$\varepsilon\in\overline{\mathbb{D}}$. Moreover, $f$ is H\"{o}lder continuous of exponent $\frac{1-\lambda}{2}$
in $\mathbb{D}$.
\ethm

\section{Coefficient estimates for the class $\mathcal{B}_H(\lambda)$}\label{LS3-sec7}

Throughout the section we consider $f=h+\overline{g}\in\mathcal{B}_{H}$, where
$$h(z)=\sum_{n=1}^{\infty}a_nz^n ~\mbox{ and }~ g(z)=\sum_{n=1}^{\infty}b_nz^n
$$
with $a_1=1$ and $\mathcal{B}_{H}$ is defined in Section \ref{LS3-sec1.3}. For $\varepsilon\in\overline{\mathbb{D}}$,
we now introduce $f_\varepsilon$ by
$$f_\varepsilon (z):=\frac{ h(z)+\varepsilon g(z)}{1+\varepsilon b_1}=\frac{1}{1+\varepsilon b_1}\sum_{n=1}^{\infty}(a_n+\varepsilon b_n)z^n.
$$
%By (\eqref{LS3-equ2.1}), below we will use the following {\it basic ~fact} repeatedly:
%$$
%f\in\mathcal{B}_{\mathcal{H}}(\lambda)\Rightarrow f_\varepsilon
%:=\frac{ h+\varepsilon g}{1+\varepsilon b_1}\in\mathcal{B}_{\mathcal{A}}\left(\frac{\lambda+1}{2}\right)~
%(\forall\varepsilon\in\overline{\mathbb{D}}).
%$$
%Note that $f_\varepsilon(z)=1/(1+\varepsilon b_1)\sum_{n=1}^{\infty}(a_n+\varepsilon b_n)z^n$.
We first determine the estimate for $a_2$.

\bthm\label{LS3-thm7.1}
If $f\in\mathcal{B}_H(\lambda)$, then we have
\be\label{LS3-equ7.1}
|a_2|\leq\frac{1}{2}\min\left\{(1-|b_1|^2)\lambda+2|b_1b_2|,
~\min_{\varepsilon\in\overline{\mathbb{D}}}\{|1+\varepsilon b_1|(\lambda+1)+2|\varepsilon b_2|\}\right\}.
\ee
%Thus, we get $|a_2|\leq\frac{\lambda+1}{2}$.
If $f\in\mathcal{B}_H^0(\lambda)$, then $|a_2|\leq \lambda /2$ and the estimate is sharp for all $\lambda>0$.
\ethm
\bpf
Let $\varepsilon\in\overline{\mathbb{D}}$ and  $f\in {\mathcal B}_H(\lambda)$ for some  $\lambda>0$. Then
$f_\varepsilon$ defined above belongs to $\mathcal{B}_{A}((\lambda+1)/2)$, and thus,
we have $|P_{f_\varepsilon}(0)|\leq||P_{f_\varepsilon}||\leq\lambda+1$ so that
$$\left |\frac{h''(0)+\varepsilon g''(0)}{h'(0)+\varepsilon g'(0)}\right |
=\left |\frac{2a_2+2\varepsilon b_2}{1+\varepsilon b_1}\right |\leq \lambda+1,
$$
which implies that $|a_2|\leq\frac{1}{2}(|1+\varepsilon b_1|(\lambda+1)+2|\varepsilon b_2|)$.
%Especially, we have $|a_2|\leq\frac{\lambda+1}{2}$ by choosing $\varepsilon=0$.
%For $\lambda\geq1$, the sharpness can be seen from the family of functions
%$F_{\frac{\lambda-1}{2},\frac{\lambda+1}{2},\theta}$ ($\theta\in[0,2\pi]$)
%since $F_{\frac{\lambda-1}{2},\frac{\lambda+1}{2},\theta}\in\mathcal{B}_{\mathcal{H}_0}(\lambda)
%\subseteq\mathcal{B}_{\mathcal{H}}(\lambda)$
%by Proposition ? and
%\begin{align*}
%H_{\frac{\lambda-1}{2},\frac{\lambda+1}{2}}=&\int_0^z\Big(\frac{1+t}{1-t}\Big)^{\frac{\lambda-1}{2}}\frac{1}{1-t}dt\\
%=&\int_0^z\Big(1+\frac{\lambda+1}{2}t+\frac{(\lambda-1)^2}{8}t^2+\cdots\Big)(1+t+t^2+\cdots)dt\\
%=&\int_0^z\Big(1+\frac{\lambda+1}{2}t+\frac{\lambda^2+2\lambda+5}{8}t^2+\cdots\Big)dt\\
%=&z+\frac{\lambda+1}{4}z^2+\frac{\lambda^2+2\lambda+5}{24}z^3+\cdots.
%\end{align*}

On the other hand, for $f\in\mathcal{B}_H(\lambda)$, it follows from \eqref{LS3-equ1.1} that
$$|P_{f}(0)|=\frac{|h''(0)\overline{h'(0)}-g''(0)\overline{g'(0)}|}{|h'(0)|^2-|g'(0)|^2}
=\frac{|2a_2-2\overline{b_1}b_2|}{1-|b_1|^2}\leq||P_{f}||=\lambda,
$$
and thus, $|a_2|\leq\frac{1}{2}((1-|b_1|^2)\lambda+2|b_1b_2|)$. Inequality  \eqref{LS3-equ7.1} follows if we combine
the last two estimates for $|a_2|$.
If $f\in\mathcal{B}_H^0(\lambda)$, then $b_1=0$ and thus, \eqref{LS3-equ7.1} reduces to  $|a_2|\leq \lambda/2$.
The function $H_{\lambda/2}=H_{\lambda/2,\lambda/2}$ defined by \eqref{LS3-equ5.2} provides the sharpness for each $\lambda>0$.
%$$H_{\lambda}(z)=\int_0^z\left (\frac{1+t}{1-t}\right )^{\lambda}dt,
%$$
%, where $H_\lambda$ is defined by \eqref{LS3-equ5.2}.
\epf

In order to indicate estimates for the coefficients of $f\in\mathcal{B}_H(\lambda)$, we consider the integral mean
$I_p(r,f)$ of $f$ defined by
%\beq \label{LS3-equ7.2}
%I_p(r):=
$$I_p(r,f)=\frac{1}{2\pi}\int_{0}^{2\pi}|f(re^{i\theta})|^p\, d\theta,
$$
where $p$ is a positive real number. Set $M_p(r,f)=(I_p(r,f))^{1/p}$,  $0<r<1$.

\bdefe\label{def7.1}
For $0<p<\infty$, the \textit{Hardy space} $H^p$ is the set of all functions $f$ analytic in $\mathbb{D}$
for which $\|f\|_p:=\sup \{M_p(r,f):\, 0<r<1\} <+\infty$, where $M_p(r,f)$ is defined as above.
%\left\{ \frac{1}{2\pi}\int_{0}^{2\pi}|f(re^{i\theta})|^p
%\,d\theta \right \}^{1/p}
%is bounded on $0<r<1$, where $I_p(r,f)$ is defined by \eqref{LS3-equ7.2}.

Let $h^p$ denote the analogous space of harmonic mappings $f$ in $\mathbb{D}$ with $\|f\|_p$ defined similarly (see \cite{dur2004}).
%The space $h^p$ consists of all harmonic mappings $f$ in $\mathbb{D}$ for which $M_p(r,f)$ $(0<r<1)$ are bounded  (see \cite{dur2004}).
\edefe

In \cite{AleMartin2015}, Aleman and Mart\'{\i}n
constructed convex harmonic mappings that do not belong to $h^{1/2}$ which settles the
question raised by Duren \cite{dur2004}. It is worth pointing out that the space $h^p$ is well-behaved for $p\geq 1$ whereas $H^p$ is comparatively
well-behaved for all $p>0$, such is not the case for $h^p$, $0<p<1$.

Since $f\in {\mathcal B}_H(\lambda)$ implies that $f_\varepsilon\in\mathcal{B}_{A}((\lambda+1)/2)$
for each $\varepsilon\in\overline{\mathbb{D}}$,  it follows from the result of \cite[p.~190]{KS2002} that
$$|a_n+\varepsilon b_n|=O\big (n^{(\lambda+1)/2-1} \big )
$$
uniformly for $\varepsilon\in\overline{\mathbb{D}}$ as $n\rightarrow \infty$ and thus, we obtain that $|a_n|+|b_n|=O\big (n^{(\lambda-1)/2}\big )$ as
$n\rightarrow \infty$. Especially, $|a_n|=O\big (n^{(\lambda-1)/2}\big)$ and $|b_n|=O\big(n^{(\lambda-1)/2}\big)$ as $n\rightarrow \infty$.
Moreover, if $\lambda<1$ and $f$ is univalent in $\mathbb{D}$, then, by Corollary \ref{LS3-cor6.1}, $f$ is  bounded and thus,
$${\rm Area}\,(f(\mathbb{D}))=\pi \Big (\sum_{n=1}^{\infty}n(|a_n|^2-|b_n|^2)\Big )<\infty,
$$
%$$\left({\rm resp. }\,{\rm Area}(h(\mathbb{D}))=\pi \Big(1+\sum_{n=2}^{\infty}n|a_n|^2\Big)<\infty,
%\quad            {\rm Area}(g(\mathbb{D}))=\pi  \sum_{n=1}^{\infty}n|b_n|^2<\infty\right),
%$$
which implies that $\sqrt{|a_n|^2-|b_n|^2}=o\big (n^{-1/2}\big )$
%$\left({\rm resp.} |a_n|=o\big (n^{-1/2}\big), |b_n|=o\big (n^{-1/2}\big)\right)$
as $n\rightarrow \infty$.

%\bee
%Let $f\in\mathcal{B}_H(\lambda)$ for some $\lambda\in[0,1]$.
%We conjecture $f$ is univalent in $\mathbb{D}$.
%\eee

%However, those estimates can be improved.
Combining the results from \cite[Section~3]{KS2002} and the implication \eqref{LS3-equ1.3}, we can get a series of results.
We omit detailed proofs, but it might be appropriate to include some necessary explanations.
In fact, we only need to modify the conditions by
replacing the parameter $\lambda$ in the theorems of \cite[Section~3]{KS2002} by $(\lambda+1)/2$ at appropriate places.

\bthm\label{LS3-thm7.2}
Let $f=h+\overline{g}\in {\mathcal B}_H(\lambda)$. Then, for any
$a>0$ and a real number $p$, we have
\be\label{LS3-equ7.2}
I_p(r, h'+\varepsilon g')=O\left ((1-r)^{-\alpha(|p|(\lambda+1)/2)-a}\right ),
\ee
for each $\varepsilon\in\overline{\mathbb{D}}$ and thus, in particular,
$$|a_n|+|b_n|=O\left (n^{\alpha((\lambda+1)/2)-1+a}\right ).
$$
For $p>0$, we get that
\be\label{LS3-equ7.3}
I_p(r, f)=O\left ((1-r)^{p-\alpha(|p|(\lambda+1)/2)-a}\right ).
\ee
Here $\alpha(\lambda)=\frac{\sqrt{1+4\lambda^2}-1}{2}$.
\ethm

%Let $M_p(r,f)=(I_p(r,f))^{1/p}$ $(p>0)$.
\bpf
The former part can be deduced from \cite[Theorem~3.1]{KS2002}.
For the later part, it follows from \eqref{LS3-equ7.2}  and \cite[Theorem~5.5]{dur1970} that
$$M_p(r, h') =O\left ((1-r)^{-(\alpha(|p|(\lambda+1)/2)+a)/p}\right )
$$
% ~\mbox{ and }~
and
$$M_p(r, h) =O\left ((1-r)^{1-(\alpha(|p|(\lambda+1)/2)+a)/p}\right ),
$$
respectively. Similar conclusions hold for $g$.
%By Theorem 5.5 in \cite{dur1970}, we know that
%$$M_p(r, h)~(M_p(r, g))=O\left ((1-r)^{1-(\alpha(|p|(\lambda+1)/2)+a)/p}\right ).
%$$
Because $M_p(r,f)\leq 2^p(M_p(r,h)+M_p(r,g))$, we finally obtain that
$$M_p(r, f)=O\left ((1-r)^{1-(\alpha(|p|(\lambda+1)/2)+a)/p}\right ),
$$
which implies \eqref{LS3-equ7.3}.
\epf

\bthm\label{LS3-thm7.3} Let
$f=h+\overline{g}\in {\mathcal B}_H(\lambda)$ for some $\lambda$ with $1.982<\lambda\leq5$.
If there exists a constant $\varepsilon\in\overline{\mathbb{D}}$ such that $h+\varepsilon g$ is univalent in $\mathbb{D}$,
then $|a_n+\varepsilon b_n|=O\big (n^{(\lambda-3)/2}\big )$ as $n\rightarrow \infty$.
In particular, if $h$ is univalent  in $\mathbb{D}$, then
$|a_n|=O\big (n^{(\lambda-3)/2}\big )$ as $n\rightarrow \infty$.
Moreover, if $h+\varepsilon g$ is univalent in $\mathbb{D}$ for  every $|\varepsilon|=1$,
then $|a_n|+|b_n|=O\big (n^{(\lambda-3)/2}\big )$ as $n\rightarrow \infty$.
The three estimates are sharp.
\ethm

\bpf
If $f\in {\mathcal B}_H(\lambda)$ for some $\lambda$ with $1.982<\lambda\leq5$, then as before we have
$f_\varepsilon\in\mathcal{B}_{A}(\frac{\lambda+1}{2})$ ($1.491<\frac{\lambda+1}{2}\leq3$)
for each $\varepsilon\in\overline{\mathbb{D}}$. The results follow from \cite[Theorem~3.2]{KS2002}.

To show the sharpness, we construct a family of functions
$$T_{\lambda,\theta}(z)=t_{\lambda}(z)+\overline{e^{i\theta}zt_{\lambda}(z)}
=\sum_{n=1}^{\infty}a_nz^n+\sum_{n=2}^{\infty}\overline{b}_n\overline{z}^n,\quad z\in\mathbb{D},
$$
where  $\lambda\in(1.982,5]$, $\theta\in\mathbb{R}$ and
$$t_{\lambda}(z)=\frac{1-(1-z)^{(1-\lambda)/2}}{(1-\lambda)/2}.
$$
First, we show that $T_{\lambda,\theta}\in{\mathcal B}_H(\lambda)$
for all $\lambda>1$ and $\theta\in\mathbb{R}$.
It suffices to prove that $||P_{T_{\lambda,\theta}}||=\lambda$ due to $T_{\lambda,\theta}\in\mathcal{H}$.
%since all of functions $T_{\lambda,\theta}$ are sense-preserving harmonic mappings with $a_1=1$.
By computation, we find that
$$P_{T_{\lambda,\theta}}(z)=\frac{1+\lambda}{2}\cdot\frac{1}{1-z}-\frac{\overline{z}}{1-|z|^2}.
$$
Also, we note that $||P_{t_{\lambda}}||=1+\lambda$. If we get $||P_{T_{\lambda,\theta}}||\leq\lambda$, then
it follows from \eqref{LS3-equ2.1} that $||P_{T_{\lambda,\theta}}||=\lambda$. Indeed we may let $z=x+iy\in\mathbb{D}$.
By computation, we obtain
\begin{align*}
&4\lambda^2|1-z|^2-\left |(1+\lambda)(1-|z|^2)-2\overline{z}(1-z)\right |^2\\
=&(\lambda-1)[(1-x)^3(1-x+\lambda(3+x))+2(3+3\lambda-2x+(1-\lambda)x^2)y^2-(\lambda-1)y^4]\\
\geq&(\lambda-1)[2(3+3\lambda-2x+(1-\lambda)x^2)y^2-2(\lambda-1)y^2]\\
\geq&2(\lambda-1)(3+3\lambda-2+1-\lambda-\lambda+1)y^2\geq0,
\end{align*}
which clearly implies that
$$(1-|z|^2)|P_{T_{\lambda,\theta}}(z)|\leq\lambda
$$
and thus, $||P_{T_{\lambda,\theta}}||\leq\lambda$. Next, we show that the functions $T_{\lambda,\theta}$ are univalent
in $\mathbb{D}$ for each $1<\lambda\leq5$ and all $\theta\in\mathbb{R}$. A simple computation shows that
$$\RE\left(1+z\frac{t''_\lambda(z)}{t'_\lambda (z)}\right)
=\RE\left(1+\frac{1+\lambda}{2}\cdot\frac{z}{1-z}\right)>1-\frac{1+\lambda}{4}\geq-\frac{1}{2},\quad z\in\mathbb{D},
$$
for $1<\lambda\leq5$. According to a well-known result, the function $t_\lambda$ is univalent and
convex in one direction (and hence, close-to-convex) in $\ID$.
Note that the dilatation of $T_{\lambda,\theta}$ is $e^{i\theta}z$.
It follows from \cite[Theorem~1]{BP2014} that the functions $T_{\lambda,\theta}$
are univalent in $\mathbb{D}$ for each $1<\lambda\leq5$ and all $\theta\in\mathbb{R}$.
Therefore, $T_{\lambda,\theta}$ is SHU and thus, $t_\lambda+\varepsilon zt_\lambda$
is univalent in $\mathbb{D}$ for each $\varepsilon\in\overline{\mathbb{D}}$.
Finally, by Stirling's formula, we have
$$|a_n|=\frac{2\Gamma((\lambda+2n-1)/2)}{(\lambda-1)n!\Gamma((\lambda-1)/2)}
\sim \frac{2}{\lambda-1}n^{(\lambda-3)/2}
~\mbox{ as $n\rightarrow\infty$}.
$$
Note that $b_n=e^{i\theta}a_{n-1}$ for each $n>1$.
Hence, $|a_n|+|b_n|=|a_n|+|a_{n-1}|=O\big (n^{(\lambda-3)/2}\big )$ as $n\rightarrow \infty$.
\epf

Given a harmonic mapping $f\in\mathcal{H}$, let $\gamma(f)$  denote the
infimum of exponents $\gamma$ such that $|a_n|+|b_n|=O\big
(n^{\gamma-1}\big )$  as $n\rightarrow \infty$, that is,
$$\gamma(f)=\overline{\lim_{n\rightarrow\infty}}\frac{\log n (|a_n|+|b_n|)}{\log n}.
$$
For the subset $X$ of ${\mathcal H}$, we let $\gamma(X)=\sup_{f\in X}\gamma(f)$.
There are some investigations about $\gamma(f)$ (resp. $\gamma(X)$) if $f$
(resp. $X$) is restricted to be analytic or special families of analytic functions.
The reader can refer to \cite{CJ, KS2002, MP} and \cite[Chapter~10]{pom1992} for some details
on this problem.

For each $\lambda\in (0,\infty)$ and $\varepsilon\in\overline{\mathbb{D}}$, we introduce
$$\mathcal{A}_{H}(\lambda,\varepsilon)=:\left \{ f_{\varepsilon}:\,  f_{\varepsilon}(z)=\frac{h(z)+\varepsilon g(z)}{1+\varepsilon g'(0)}
~\mbox{ and }~f=h+\overline{g}\in{\mathcal B}_H(\lambda)\right\}
$$
and obtain the following theorem.

\bthm\label{LS3-thm7.4}
For each $\lambda\in (0,\infty)$ and $\varepsilon\in\overline{\mathbb{D}}$, we have
$$\max\{(\lambda-1)/2,0\}\leq\gamma(\mathcal{B}_H(\lambda))\quad
\left(\gamma(\mathcal{A}_{H}(\lambda,\varepsilon))\right)\leq \alpha((\lambda+1)/2),
$$
where $\alpha(\lambda)=\frac{\sqrt{1+4\lambda^2}-1}{2}$. In particular,
$\gamma({\mathcal B}_H(\lambda))=O((\lambda+1)^2)$ and
$\gamma(\mathcal{A}_{H}(\lambda,\varepsilon))=O((\lambda+1)^2)$
as $\lambda\rightarrow 0$.
\ethm

We continue the discussion by mentioning a connection with integral means for univalent analytic functions.
For a univalent harmonic mapping $f=h+\overline{g}\in {\mathcal S}_H$,
a complex number $\varepsilon\in\overline{\mathbb{D}}$ and a real number $p$, we let
$$\beta_{f_\varepsilon}(p)=\overline{\lim_{r\rightarrow 1^{-}}}\frac{\log
I_p(r, f'_\varepsilon)}{\log \frac{1}{1-r}}.
$$
%where $f_{\varepsilon}=\frac{h+\varepsilon g}{1+\varepsilon g'(0)}$.
Clearly, for a univalent analytic function $f\in\mathcal{A}\cap\mathcal{S}_{H}$,
$$\beta_{f}(p)=\overline{\lim_{r\rightarrow 1^{-}}}
\frac{\log I_p(r, f')}{\log \frac{1}{1-r}}.
$$
Brennan conjectured that $\beta_f(-2)\leq 1$ for univalent analytic functions $f$ (see \cite[Charpter~8]{pom1992}).

As a corollary to Theorem \ref{LS3-thm7.2}, we have

\bthm\label{LS3-thm7.5}
For $f\in\mathcal{B}_{H}(\lambda)$ and a real number $p$,
$$\beta_{f_\varepsilon}(p)\leq \alpha(|p|(\lambda+1)/2)=\frac{\sqrt{1+p^2(1+\lambda)^2}-1}{2}
$$
holds for each  $\varepsilon\in\overline{\mathbb{D}}$.
In particular, the Brennan conjecture is true for  every univalent harmonic mapping $f$ with $\|P_f\|\leq\sqrt{2}-1$.
\ethm

\section{The space ${\mathcal B}_H(\lambda)$ and the Hardy space}\label{LS3-sec8}

For a harmonic mapping $f=h+\overline{g}$ in $\mathbb{D}$, the
Bloch seminorm is given by (see \cite{col1989-1})
$$\|f\|_{\mathscr{B}_H}=\sup_{z\in \mathbb{D}} (1-|z|^2)\big(|h'(z)|+|g'(z)|\big),
$$
and $f$ is called a (harmonic) Bloch mapping when $\|f\|_{\mathscr{B}_H}<\infty$.
%In the recent years, the class of harmonic Bloch mappings has been studied extensively together with its
%higher dimensional analog (see for example, \cite{CPVW,CPW2011,CPW2012,col1989-1,col1989-2} and the references therein).
Let ${\rm BMOA}$ (resp. ${\rm BMOH}$) denote the class of analytic functions
(resp. harmonic mappings) that have bounded mean oscillation on the
unit disk $\mathbb{D}$ (see \cite{abu}).
Kim \cite{kim} showed some relationships among $\mathcal{B}_{A}(\lambda)$, $H^p$ and BMOA (see also \cite{KS2009}).
Combined with the study on Bloch, BMO and  univalent harmonic mappings (see \cite{abu}),
a generalization of Kim's result is given in \cite{PQW}. Basic properties about analytic Bloch functions may be
obtained from \cite{ACP,pom1992}.

Our results are based on the following observation.
It follows from Theorem \ref{LS3-thm6.1} (6) that the inequality
$$ |f(z)|\leq (1+|b_1|)\int_{0}^{|z|}\left (\frac{1+t}{1-t}\right )^{\frac{\lambda+1}{2}}\,dt, \quad z\in\mathbb{D},
$$
holds for every $f\in\mathcal{B}_H(\lambda)$, which implies that

\bei
\item $f$ is bounded when $\lambda<1$,

\item $f(z)=O(-\log(1-|z|))$ $(|z|\rightarrow 1)$ when $\lambda=1$,
and

\item $f(z)=O((1-|z|)^{1-\frac{\lambda+1}{2}})$ $(|z|\rightarrow 1)$ when
$\lambda>1$.
\eei

%Associated with the above results and \eqref{LS3-equ1.3},
%we acquire the following theorems by the similar proofs of theorems in Section 4 of \cite{PQW}.
On the other hand, the proofs of our results are similar to that of results of \cite[Section~4]{PQW}.
%Let's recall partial works in \cite{PQW}.
Let
$${\mathcal{T}}_{H}(\lambda)=\{f=h+\overline{g}\in {\mathcal{H}}:\, \|T_f\|\leq 2\lambda\}
$$
 with
\beqq
\|T_f\|:=\sup_{z\in \mathbb{D},\ \theta\in [0,
2\pi]} (1-|z|^2)\left
|\frac{h''(z)+e^{i\theta}g''(z)}{h'(z)+e^{i\theta}g'(z)}\right|.
\eeqq
For the one parameter family ${\mathcal{T}}_{H}(\lambda)$, the authors %not only
%obtained sharp distortion, growth and covering theorems, and the growth of coefficients,
showed its relationship with Hardy spaces in \cite[Section~4]{PQW}. %, which is a generalization of \cite{kim}.
%Thus, they extended most of results in \cite{KS2002} to harmonic mappings.
%However, the norm attached to ${\mathcal{T}}_{H}(\lambda)$ is illegibility for a given harmonic mapping.
%Therefore, their results have some limitations.
Note that $||T_f||=\sup_{\theta\in [0,2\pi]} ||P_{h+e^{i\theta}g}||.$
If $f\in\mathcal{B}_H(\lambda)$, then it is easy to see  that $f\in\mathcal{T}_H(\frac{\lambda+1}{2})$ from \eqref{LS3-equ2.1}.
%Therefore, Theorems \ref{LS3-thm8.1} and \ref{LS3-thm8.2} below follows from the proofs of theorems in \cite[Section4]{PQW}.
Therefore, applying the above observation and replacing $h+e^{i\theta}g~(\theta\in[0,2\pi])$ (resp. $\lambda$)
%in Lemma 4.1 (resp. the proof of (3) in Theorem 4.1 and Theorem 4.2) in \cite{PQW}
to $h+\varepsilon g~(\varepsilon\in\overline{\mathbb{D}})$ (resp. $(\lambda+1)/2$) in corresponding proofs of \cite[Section~4]{PQW},
then we can easily obtain the following results. So we omit their proofs.

\bthm\label{LS3-thm8.1}
\begin{enumerate}
\item If $\lambda<1$, then $\mathcal{B}_H(\lambda)\cap {\mathcal{S}}_H\subset
h^{\infty}$.

\item If $\lambda=1$, then $\mathcal{B}_H(\lambda)\cap {\mathcal{S}}_H\subset
{\rm BMOH}$.

\item If $\lambda>1$, then $\mathcal{B}_H(\lambda)\cap {\mathcal{S}}_{H_K} \subset
h^{p}$ for every $0<p<2/(\lambda-1)$, where $K\geq1$ and
$${\mathcal{S}}_{H_K} =\{f=h+\overline{g}\in {\mathcal{S}}_H:\, f  ~\mbox{ is $K$-quasiconformal}\}.
$$
\end{enumerate}
\ethm

\bthm\label{LS3-thm8.2}
Let $\lambda\geq 1$. Then $\mathcal{B}_H(\lambda)\subset h^p$ with
$0<p<p_0(\lambda )=\frac{4}{(\lambda +3)(\lambda-1)}$,
where $ p_0(\lambda )=\infty$ if $\lambda=1$.
\ethm

\br
Theorem \ref{LS3-thm8.2} can be directly obtained by choosing
$p-\alpha(p(\lambda+1)/2)>0$ in Theorem \ref{LS3-thm7.2}.
\er

\bcor \label{LS3-cor8.1}
A uniformly locally univalent harmonic mapping $f$ in
$\mathbb{D}$ is contained in the Hardy space $h^p$ for some
$p=p(f)>0$.
\ecor

%In \cite{Ki}, Kim also conjectured that the assertion (3) in Theorem
%\Ref{Thm a} holds for ${\mathcal{B}}(\lambda)$.

\section{Subordination principles for the estimate of PSN} \label{LS3-sec9}

In this section, $\mathcal{A}_D$ denotes the class of analytic functions $\phi$ from $\mathbb{D}$ into itself
and  $\mathcal{A}_D^0$ denotes the subclass of $\mathcal{A}_D$ with the normalization  $\phi(0)=0$.
%First we recall some definitions of subordination for analytic functions defined in the unit disk.
If $f$ and $F$ are restricted to be analytic, then we say that  $f$ is said to be
 {\it subordinate} (resp. {\it weakly subordinate}) to $F$
(written $f\prec F$ (resp. $f\preceq F$))
if there exists  a function $\phi\in\mathcal{A}_D^0$ (resp. $\phi\in\mathcal{A}_D$) such that $f(z)=F(\phi(z))$ in $\mathbb{D}$.

In 2000, Schaubroeck \cite{sch} generalized the notion of subordination to harmonic mappings.
A harmonic mapping $f$ is subordinate to a harmonic mapping $F$, still denoted
by $f\prec F$, if there is a function $\phi\in\mathcal{A}_D^0$ such that $f= F\circ \phi$.
%Suppose $f = h+\overline{g}$ and $F = H+\overline{G}$,
%where $h, ~g, ~H$ and $G$ are analytic in $\mathbb{D} $ with $h(0) = H(0)$ and $g(0) = G(0)$.
%If $f\prec F$, then it is easy to see that $h \prec H$ and $g \prec G$.

Note that if the analytic function $F$  is univalent in $\mathbb{D}$,  then  $f\prec F$ if and only if
that $f(0)=F(0)$ and $f(\mathbb{D})\subseteq F(\mathbb{D})$.
However, this property is not true for harmonic mappings.
As in \cite{mui},  a harmonic mapping $f$ is said to be weakly subordinate to the harmonic mapping $F$
if $f(\mathbb{D})\subseteq F(\mathbb{D})$.
%Particularly, if $f$ and $F$ are chosen by analytic functions and $f$ is weakly subordinate to  $F$, which equivalents to  that
%$f$ is subordinate to $F$.
%Such notation is not generalization of weakly subordination for analytic functions in \cite{KS2002}.

In this article, $f=h+\overline{g}\preceq F=H+\overline{G}$
means that there exists  a function $\phi\in\mathcal{A}_D$ such that $h=H\circ \phi$ and $g=G\circ \phi$.
Clearly, if $f\preceq F$, then $f$ is weakly subordinate to $F$ in the sense of Muir.
The following result is a generalization of   \cite[Theorem~4.1]{KS2002}.

\bthm\label{LS3-thm9.1}
{\rm (Subordination principle I)} Let $f=h+\overline{g}$  be a harmonic mapping in $\mathbb{D}$
and $F=H+\overline{G}\in\mathcal{B}_{H}$.
If $h'+\overline{g'}\preceq H'+\overline{G'}$, then we have $||P_f||\leq||P_F||$.
In this case, $f$ is {\rm ULU} in $\mathbb{D}$.
\ethm
\bpf
By assumption, there exists a function $\phi\in\mathcal{A}_D$ such that  %$f= F\circ \phi$. Thus,
%$$f_z(z)= H'(\phi(z))\phi'(z) ~\mbox{ and }~f_{\overline{z}}(z)= \overline{G'(\phi(z))\phi'(z)}.
%$$
$h'=H'\circ \phi$ and $g'=G'\circ \phi$.
Therefore, $f$ is sense-preserving in $\mathbb{D}$ since $F$ is sense-preserving in $\mathbb{D}$.
Moreover, we have $P_h=(P_H\circ \phi)\phi'$  and
$$\omega_f(z)=\frac{g'(z)}{h'(z)}=\frac{G'(w)}{H'(w)}=\omega_F(w),\quad w= \phi(z).
$$
Consequently,
$$\frac{\overline{\omega_f(z)}\omega'_f(z)}{1-|\omega_f(z)|^2}
=\frac{\overline{\omega_F(w)}\omega'_F(w)}{1-|\omega_F(w)|^2}\phi'(z).
$$
It follows from (\ref{LS3-equ1.1}) that $P_f=(P_F\circ \phi)\phi'$. By Schwarz-Pick's lemma,
$$|\phi'(z)|\leq\frac{1-|\phi(z)|^2}{1-|z|^2}
$$
and using this, we find that
$$(1-|z|^2)|P_f(z)|=(1-|z|^2)|\phi'(z)P_F(\phi(z))|\leq(1-|\phi(z)|^2)|P_F(\phi(z))|\leq||P_F||.
$$
The desired conclusion follows.
\epf

%\br
%Theorem \ref{LS3-thm9.1} continues to hold if the assumption $h'+\overline{g'}\preceq H'+\overline{G'}$ is replaced by
%$f=h+\overline{g}\preceq F=H+\overline{G}$. This is because the last condition implies that $f= F\circ \phi$ so that
%$$f_z(z)= H'(\phi(z))\phi'(z) ~\mbox{ and }~f_{\overline{z}}(z)= \overline{G'(\phi(z))\phi'(z)}.
%$$
%\er
%\br
%If the condition $h'+\overline{g'}\preceq H'+\overline{G'}$
%in Theorem \ref{LS3-thm9.1} is replaced by $h'(0)=H'(0)$,
%$g'(0)=G'(0)$ and $h'+\overline{g'}\prec H'+\overline{G'}$,
%we still have $||P_f||\leq||P_F||$.
%\er

Often, the property of a sense-preserving harmonic mapping
is mainly decided by its analytic part.
As another example of it, we have
%Therefore, by means of subordination,  to acquire the relationships
%among the pre-Schwarzian norms of sense-preserving harmonic mappings,
%it is more easily to consider that of their analytic parts.

\bthm\label{LS3-thm9.2}
{\rm (Subordination principle II)}
Let $f=h+\overline{g}$ be a sense-preserving harmonic mapping in $\mathbb{D}$
and $F=H+\overline{G}\in\mathcal{B}_{H}$ such that $h'\preceq H'$.
Then we have $||P_f||\leq||P_F||+2$.
Thus, $f$ is {\rm ULU} in $\mathbb{D}$.
\ethm
\bpf
Since $F\in\mathcal{B}_H$,
we know that $H\in\mathcal{B}_{A}$ by Theorem \ref{LS3-thm4.1}.
Clearly, $\mathcal{B}_{A}\subseteq\mathcal{B}_H$.
It follows from the assumption and Theorem \ref{LS3-thm9.1} that $||P_h||\leq||P_H||.$
Using the inequality (\ref{LS3-equ2.1}) twice, we obtain that
$$||P_f||\leq||P_h||+1\leq||P_H||+1\leq||P_F||+2
$$
and the proof is complete.
\epf

Similar to Theorem \ref{LS3-thm9.2}, few other results on subordination of analytic functions in \cite[Section~4]{KS2002}
can be transplanted to the case of sense-preserving harmonic mappings
by considering its analytic parts.% of harmonic mappings and combine the Proposition \ref{LS3-prop1} appropriately.
%For example, if $f=h+\overline{g}\in\mathcal{B}_{\mathcal{H}}$ satisfies $\RE(h'+\varepsilon g')>0$ for some
%We leave   it to the reader as an excise.
%If the condition $h'\preceq H'$ in Corollary \ref{LS3-cor6} is changed to $h\preceq H$,
%then there is a $\phi\in\mathcal{AD}$ such that $h=H(\phi)$. By computation, we have
%$$P_h=\frac{h''}{h'}=\frac{H''(\phi)(\phi')^2+H'(\phi)\phi''}{H'(\phi)\phi'}=P_H(\phi)\phi'+P_\phi.$$

\subsection*{Acknowledgments}
The work was completed during the visit of the first author to the Indian Statistical Institute,
Chennai Centre and this author thanks the institute for the support and the hospitality.
The research of the first author was supported by the NSFs of China (No. 11571049),
the Construct Program of the Key Discipline in Hunan Province,
the Science and Technology Plan Project of Hunan Province (No. 2016TP1020) and
the Science and Technology Plan Project of Hengyang City (2017KJ183). The  work of the second author is
supported by Mathematical Research Impact Centric Support of DST, India  (MTR/2017/000367).
The second author is currently at Indian Statistical Institute (ISI), Chennai Centre, Chennai, India.
The authors thank the referee for his/her comments.

%\subsection*{Conflict of Interests}
%The authors declare that there is no conflict of interests regarding the publication of this paper.

\end{document}